# Computational multiphase periporomechanics for unguided cracking in unsaturated porous media

Shashank Menon, Xiaoyu Song[1]

*Engineering School of Sustainable Infrastructure & Environment University of Florida, Gainesville, FL 32611, USA.*

**Abstract**

In this article we formulate and implement a computational multiphase periporomechanics model for unguided fracturing in unsaturated porous media. The same governing equation for the solid phase applies on and off cracks. Crack formation in this framework is autonomous, requiring no prior estimates of crack topology. As a new contribution, an energy-based criterion for arbitrary crack formation is formulated using the peridynamic effective force state for unsaturated porous media. Unsaturated fluid flow in the fracture space is modeled in a simplified way in line with the nonlocal formulation of unsaturated fluid flow in the bulk. The formulated unsaturated fracturing periporomechanics is numerically implemented through an implicit fractional step algorithm in time and a two-phase mixed meshless method in space. The two-stage operator split converts the coupled periporomechanics problem into an undrained deformation and fracture problem and an unsaturated fluid flow in the deformed skeleton configuration. Numerical simulations of in-plane open and shear cracking are conducted to validate the accuracy and robustness of the fracturing unsaturated periporomechanics model. Then numerical examples of wing cracking and nonplanar cracking in unsaturated soil specimens are presented to demonstrate the efficacy of the proposed multiphase periporomechanics model for unguided cracking in unsaturated porous media.

*Keywords:* unsaturated porous media, unguided cracking, fracture fluid flow, nonlocal, periporomechanics

## 1. Introduction

The mechanical and physical behavior of unsaturated porous media (e.g., unsaturated soils) play a significant role in resilience and sustainability of civil infrastructures (e.g., [1–4]). Cracking in unsaturated soils can significantly deteriorate and compromise the integrity of civil infrastructures on such materials ([5– 8]). For instance, the volume shrinkage by variations of matric suction (i.e., the difference between pore air and water pressures) could generate tensile cracks in unsaturated soils. The bearing capacity of structural foundations on such soils can be significantly reduced by the presence of tensile cracks. Surface cracks caused by rapid drawdown of a reservoir can be a trigger for landslides [5]. Desiccation cracking in clay can dramatically increases its hydraulic conductivity that will compromise its ability to act as liner material for landfills [9]. With advance of supercomputers, computational methods play an increasingly significant role in modeling the mechanical and physical behavior including cracking of unsaturated soils (e.g., [4]). In this study, different from the celebrated computational method such as the extended finite element method (e.g., [10, 11]) and as a new contribution, we develop a fully coupled nonlocal mathematical framework for modeling fracturing and fluid flow in unsaturated porous media using periporomechanics (e.g., [12–14] and see Section 2.1 for a brief review) that is a nonlocal reformulation of classical poromechanics [15–17] through

---

[1] Corresponding author

*Email address:* xysong@ufl.edu (Xiaoyu Song )

peridynamics [18, 19] for modeling variably saturated porous media. We refer to the celebrated literature (see [6, 7], among many others) on other robust computational techniques on modeling cracking in fracturing porous media under variably saturated conditions. Note that the original peridynamics was proposed by Silling for modeling deformation and cracks in single-phase solid materials. In the proposed periporomechanics model the same momentum balance equation applies on and off the crack surface by using the effective force state concept [12], as in the classical peridynamics for modeling cracks in solids. Next, we briefly review the application of peridynamics for modeling cracks focusing on the criterion of initiation and propagation of cracks.

In peridynamics, cracks nucleate, grow, branch, and merge when and where it is energetically favorable for them to do so according to the global equations and constitutive model ([20, 21]). This salient feature of autonomous crack growth obviates the need for the remedial techniques of classical fracture mechanics. We refer to the literature (e.g., [22]) for an in-depth comparison of peridynamics and classical fracture mechanics [23]. Peridynamics has been applied to model fractures in solids (e.g., [24–30], among many others) and porous media (e.g., [8, 31–35], among others). There are two criteria for the inception and propagation of cracks, namely "bond-breakage" using either a kinematic criterion (stretch, shear) (e.g., [8, 24, 36–39]) and a strain energy criterion (e.g., [33, 40–43]). In the kinematic-based criterion, peridynamic bonds are assumed to fail irreversibly when the relative extension exceeded a pre-determined limit value called "critical stretch". This criterion has some limitations including its ad-hoc feature, its inability to correctly model shear cracks, and its apparent lack of any rational basis or theoretical justification for the link between critical stretch and crack energy release rate. In the energy-based criterion, the elastic strain energy is correlated to the fracture potential energy. Madenci and Oterkus [42] proposed an energy-based failure criterion for fracture, based on the peridynamic equivalent of the classical $J$-integral [44]. Breitenfeld et al. [41] formulated a peridynamic equivalent of the $J$-integral to extract the classical stress intensity factors of linear elastic fracture mechanics as a means of examining the stress singularities. Lipton et al. [43] developed an elastic peridynamic material based on a novel strain energy density function capable of modeling cracks. In this method, cracks are modeled as material instability in the softening regime allowing spontaneous nucleation of cracks without the assistance of supplemental criteria. Note that both criteria have been utilized to model fracture in porous media (e.g., [8, 34, 35, 45]) through the effective stress concept for porous media [1, 15]. In this article, as a new contribution, in the fracturing unsaturated periporomechanics model we developed an energy-based crack criterion through the recently proposed effective force state concept [12] for unsaturated porous media.

In [14], the authors developed a stable coupled periporomechanics model for dynamic analysis of saturated porous media with no cracks. In the present article, as a new contribution the formulation in [14] is extended to model unguided fracturing and fluid flow in unsaturated porous media. In this fracturing unsaturated periporomechanics model, the same governing equation for the solid phase applies on and off cracks. In line with the original peridynamics [18, 19, 46], crack formation in this coupled framework is autonomous, requiring no prior estimates of crack topology. As a new contribution, the criterion for arbitrary crack formation is formulated based on the peridynamic effective force state for unsaturated porous media. Unsaturated fluid flow in the fracture space is coupled to the fluid flow in bulk through a leak-off term. Different from the fully coupled numerical implementation in [14], the formulated unsaturated fracturing periporomechanics is numerically implemented through the celebrated fractional step algorithm (also called staggered algorithm in the literature, see [6, 15, 16, 34, 45, 47–49] and others) for computational efficiency and accuracy. The two-stage operator split converts the coupled fracturing periporomechanics problem into an undrained deformation and fracture problem and an unsaturated fluid flow in the deformed configuration. Note that we refer to the literature for explicit or implicit implementation of peridynamics for modeling cracks in single-phase solids (e.g., [41, 46, 50–53], and among others). Numerical simulations are conducted to validate the accuracy and robustness of the computational fracturing unsaturated periporomechanics model. We run numerical simulations of crack branching in porous media and perform sensitivity analysis of the results with respect to hydraulic conductivity and initial matric suction. Finally, we demonstrate the efficacy of the proposed fracturing periporomechanics model for modeling non-planar cracks by simulating non-planar cracks triggered by varying matric suction in a three-dimensional soil specimen.

The contribution of this article includes (i) the formulation of the energy-based crack criterion based on the effective force state concept for unsaturated porous media, (ii) the implicit fractional-step and mixed meshless numerical implementation of the fracturing unsaturated periporomechanics model in time and space, and (iii) the validation of the numerical implementation by comparing numerical results from this newly formulated method with the classical extended finite element method and demonstration of the new coupled framework for modeling non-planar cracks in unsaturated soils. For sign convention, the assumption in continuum mechanics is followed, i.e., for solid skeleton tension is positive and compression is negative, and for pore fluid compression is positive and tension is negative.

## 2. Unsaturated fracture periporomechanics for unguided cracking

### 2.1. Unsaturated periporomechanics

Periporomechanics is a reformulation of classical poromechanics through peridynamic state concept [19] for modeling continuous or discontinuous deformation and fluid flow in variably saturated porous media ([12–14]). In periporomechanics, it is hypothesized that a porous material body can be represented by a finite number of mixed material points that are endorsed with two kinds of degree of freedom, i.e., solid displacement and fluid pressure. It is a fully coupled strong nonlocal theory in that a material point $x$ has direct poromechanical and physical interactions with any material point $x'$ in its nonlocal *family* H [19], a spherical domain centered at $x$. The radius of H denoted by $\delta$ is called the *horizon* of the porous media. For conciseness of notations, it is assumed that a peridynamic state variable without a prime is evaluated at $x$ on the associated bond $\xi = x' - x$ and the peridynamic state variable with a prime is evaluated at $x'$ on the associated bond $\xi' = x - x'$, e.g., $\underline{T} = \underline{T}[x]<x'-x>$ and $\underline{T}' = \underline{T}[x]<x-x'>$. Similarly, a non-state variable without a prime is associated with $x$ and that with a prime is associated with $x'$, e.g., $S_r$ is the degree of saturation at $x$ and $S_r'$ is the degree of saturation at $x'$.

In unsaturated periporomechanics [12], under the assumption of passive air pressure the balance of linear momentum at $x$ reads

$$\int_{\mathcal{H}} \left[ \left( \underline{T} - S_r \underline{T}_w \right) - \left( \underline{T}' - S_r' \underline{T}_w' \right) \right] \, dV' + \rho \boldsymbol{g} = \rho \ddot{\boldsymbol{u}} \tag{1}$$

where $\ddot{\boldsymbol{u}}$ is the acceleration, $\underline{T}$ is the force vector state at $x$ along a bond $\xi$, $\underline{T}$ is the effective force vector state at $x$ along a bond $\xi$, $\underline{T}_w$ is the pressure force vector state at $x$ along a bond $\xi$, $S_r$ is the degree of saturation of water, $\mathcal{H}$ is the family or the neighborhood of the point at $x$, $\boldsymbol{g}$ is the acceleration due to gravity, $dV'$ is the volume of the neighbor $x'$ and $\rho$ is the density of unsaturated porous media,

$$\rho = \rho_s(1 - \phi) + S_r \phi \rho_w \tag{2}$$

where $\rho_s$ is solid density, $\rho_w$ is water density and $\phi$ is porosity. Under the assumptions of incompressible solid grain and water the balance of mass at $x$ reads

$$\phi \frac{dS_r}{dt} + S_r \int_{\mathcal{H}} (\underline{\dot{\mathcal{V}}}_s - \underline{\dot{\mathcal{V}}}_s') \, dV' + \frac{1}{\rho_w} \int_{\mathcal{H}} (Q - Q') \, dV' + Q_s = 0, \tag{3}$$

where $\dot{\mathcal{V}}_s$ and $\dot{\mathcal{V}}_s'$ are the solid volume change rate states at $x$ and $x'$, respectively, $Q$ and $Q'$ are the fluid flow states at $x$ and $x'$, respectively, and $Q_s$ is a source/sink term.

The deformation state and pore fluid pressure state are essential state variables to construct constitutive models in periporomechanics [13]. The deformation state $\underline{Y}$ at $x$ reads

$$\underline{Y} = \boldsymbol{y}' - \boldsymbol{y} \tag{4}$$

where $\boldsymbol{y}'$ and $\boldsymbol{y}$ are the deformed positions of $x'$ and $x$ respectively, respectively. Let $\boldsymbol{u}'$ and $\boldsymbol{u}$ be displacements at $x'$ and $x$, respectively,

$$\boldsymbol{y} = \boldsymbol{x} + \boldsymbol{u}, \quad \boldsymbol{y}' = \boldsymbol{x}' + \boldsymbol{u}'. \tag{5}$$

The fluid flow/pressure potential state at $x$ reads

$$\underline{\Phi} = p'_w - p_w, \tag{6}$$

where $p_w$ and $p'_w$ are water pressures at $x$ and $x'$, respectively.

In line with the stabilized multiphase correspondence principle [14], the effective force state and unsaturated fluid flow state can be written as follows. First, the effective force state with stabilization reads

$$\underline{T} = \omega \left( \overline{P} K^{-1} \xi + \frac{GC}{\omega_0} \underline{\mathcal{R}}^s \right), \quad C = \frac{18K}{\pi \delta^4}, \tag{7}$$

where $\omega$ is the influence function (also called weighting function), $\overline{P}$ is the effective Piola-Kichhoff stress tensor, $K$ is the shape tensor [19], $\underline{\mathcal{R}}^s$ is the residual deformation state, and $G$ and $K$ are the shear and bulk moduli of the skeleton, and $\omega_0$ is defined as

$$\omega_0 = \int_{\mathcal{H}} \omega \, dV'. \tag{8}$$

The effective Piola-Kirchhoff stress tensor [54, 55] is defined as

$$\overline{P} = P + J(p_w \mathbf{1}) \widetilde{F}^{-T}, \tag{9}$$

where $\overline{P}$ is the total Piola-Kirchhoff stress tensor, $\mathbf{1}$ is the second-order identity tensor, and $J$ is the determinant of $\widetilde{F}$ that is the nonlocal deformation gradient

$$\widetilde{F} = \left( \int_{\mathcal{H}} \omega \, \underline{Y} \otimes \xi \, dV' \right) K^{-1}. \tag{10}$$

The residual deformation state $\underline{\mathcal{R}}^s$ is defined as

$$\underline{\mathcal{R}}^s = \underline{Y} - \widetilde{F} \xi. \tag{11}$$

It is noted that for a uniform deformation the residual deformation state $\underline{\mathcal{R}}^s$ is null. Given $\widetilde{F}$ the effective Piola-Kichhoff stress tensor $\overline{P}$ can be determined through the classical constitutive models for unsaturated soils. Let $\phi_0$ be the initial porosity, the porosity in (3) can be written as [54]

$$\phi = 1 - (1 - \phi_0)/J. \tag{12}$$

Similarly, the fluid flow state with stabilization at $x$ [12] is written as

$$\underline{Q} = \omega \left( \rho_w q K^{-1} \xi + \frac{GK_p}{\omega_0} \underline{\mathcal{R}}^w \right), \quad K_p = \frac{6\rho_w k_w}{\pi \delta^4}. \tag{13}$$

where $q$ is the fluid flux vector, $\underline{\mathcal{R}}^w$ is the residual pressure potential state, $K_p$ is the fluid stabilization parameter, $k_w$ is the intrinsic permeability, and $\mu_w$ is water viscosity. The fluid flux $q$ can be determined by the generalized Darcy's law for unsaturated fluid flow as

$$q = -k^r k_w \widetilde{\nabla \Phi}, \quad \widetilde{\nabla \Phi} = \left( \int_{\mathcal{H}} \omega \, \underline{\Phi} \, \xi \, dV' \right) K^{-1}, \tag{14}$$

where $k^r$ is the relative permeability and $\widetilde{\nabla \Phi}$ is the nonlocal fluid pressure gradient. The residual fluid potential state $\underline{\mathcal{R}}^w$ is defined as

$$\underline{\mathcal{R}}^w = \underline{\Phi} - \widetilde{\nabla \Phi} \xi. \tag{15}$$

Note that $\underline{\mathcal{R}}^w$ is null for a uniform fluid pressure potential. Here, the soil-water retention curve [56] is described by the celebrated van Genuchen equation [57] as

$$S_r = S_1 + (1 - S_2) \left[1 + \left(\frac{-p_w}{s_a}\right)^n\right]^{(1-n)/n} \tag{16}$$

where $s$ is matric suction, i.e., air pressure minus water pressure (or negative water pressure assuming passive air pressure), $S_1$, $S_2$, $s_a$ and $n$ are material constants. Given $S_r$ the relative permeability $k_r$ can be written as

$$k^r = S_r^{1/2} \left[1 - (1 - S_r^{1/m})^m\right]^2 \tag{17}$$

where $m = (n - 1)/n$.

Next, we present a formulation of so-called ordinary material models for unsaturated porous media that do not exhibit the same instability observed in the original correspondence material models. Here the ordinary means that the force state is parallel to the bond $\underline{\xi}$. Let $\underline{y} = |\underline{Y}|$, the total force state reads

$$\underline{T} = (\underline{t} - S_r \underline{t_w}) \left(\frac{\underline{Y}}{\underline{y}}\right) \tag{18}$$

where $\underline{t}$ and $\underline{t_w}$ are the scalar effective force state and the scalar water pressure state, respectively. For a poroelastic material, the scalar effective force state $\underline{t}$ can be defined as

$$\underline{t} = \underline{\omega} \left(\frac{3K}{m_v} \theta \underline{x} + \frac{15G}{m_v} \underline{e}^d\right) \tag{19}$$

where $\underline{x} = |\underline{\xi}|$, $\theta$ is the dilatation, $\underline{e}^d = \underline{e} - (\theta \underline{x}/3)$ is a measure of deviatoric deformation, $\underline{e} = \underline{y} - \underline{x}$, $m_v$ is the weighted volume, and $K$ and $G$ are the elastic bulk and shear moduli respectively as introduced previously.

$$m_v = \int_{\mathcal{H}} \underline{\omega} \, \underline{x}^2 \, dV' \tag{20}$$

$$\theta = \frac{3}{m_v} \int_{\mathcal{H}} \underline{e} \, \underline{x} \, dV' \tag{21}$$

In line with equation (19), the water pressure scalar state in (18) can be defined as

$$\underline{t_w} = \underline{\omega} \frac{3 p_w \underline{x}}{m_v} \tag{22}$$

From (18) and (19), the effective force state through the ordinary material model for the solid skeleton reads

$$\overline{\underline{T}} = \underline{t} \left(\frac{\underline{Y}}{\underline{y}}\right) \tag{23}$$

We note that the ordinary material model in (23) is a special case of the general correspondence material model in (7). To support this statement, we first demonstrate that the shape tensor can be expressed by the weighted volume in what follows. The spherical coordinate system in Figure 1 is adopted to facilitate the derivation. Given the notations in Figure 1, the components of $\underline{\xi}$ in the Cartesian coordinate system read

$$\xi_1 = \underline{x} \sin \overline{\alpha}_1 \sin \overline{\alpha}_2, \tag{24}$$

$$\xi_2 = \underline{x} \cos \overline{\alpha}_1, \tag{25}$$

$$\xi_3 = \underline{x} \sin \overline{\alpha}_1 \cos \overline{\alpha}_2. \tag{26}$$

Let $\underline{\omega} = 1$ the weighted volume $m_v$ is written as

$$m_v = \int_0^\delta \int_0^{2\pi} \int_0^\pi \underline{x}^4 \sin \overline{\alpha}_1 \, d\overline{\alpha}_1 \, d\overline{\alpha}_2 \, d\underline{x}$$
$$= \frac{4\pi \delta^5}{5} \tag{27}$$

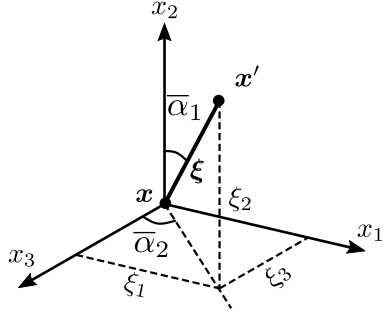

Figure 1: Two material points in the spherical coordinate system.

Similarly, the shape tensor can be written as

$$K_{ij} = \int_{\mathcal{H}} \underline{\omega} \xi_i \xi_j dV'$$
$$= \int_0^{\delta} \int_0^{2\pi} \int_0^{\pi} \underline{\omega} \, \xi_i \xi_j \underline{x}^2 \sin \overline{\alpha}_1 \, d\overline{\alpha}_1 \, d\overline{\alpha}_2 \, d\underline{x}, \quad (28)$$

where $i,j$ = 1,2,3. It follows from (24), (25), (26) and (28) that for $\underline{\omega}$ = 1 we can readily show

$$K_{ij} = \frac{4\pi\delta^5}{15} \delta_{ij}, \quad (29)$$

From (27) and (29), the shape tensor can be written as

$$\boldsymbol{K} = (m_v/3)\mathbf{1}. \quad (30)$$

Given (30), we can establish the equivalence between correspondence material models and ordinary material models for unsaturated porous media in what follows. Let us assume an isotropic deformation of unsaturated porous media,

$$\underline{\boldsymbol{Y}} = (1 + \epsilon)\boldsymbol{\xi}, \quad (31)$$

where is a constant scalar and $|\epsilon| \ll 1$. Then from (10) and (29), the nonlocal deformation gradient can be written as

$$\widetilde{\boldsymbol{F}} = (1 + \epsilon)\mathbf{1}. \quad (32)$$

It follows from (7), (9) and the assumption $|\epsilon| \ll 1$ that the total force state can be written as

$$\underline{\boldsymbol{T}} = \underline{\omega} \left( \overline{\boldsymbol{\sigma}} + S_r p_w \mathbf{1} \right) (1) \left( \frac{3}{m_v} \mathbf{1} \right) \boldsymbol{\xi}$$
$$= \underline{\omega} \left( 3K\epsilon\mathbf{1} + S_r p_w \mathbf{1} \right) \left( \frac{3}{m_v} \mathbf{1} \right) \underline{x} \frac{(1+\epsilon)\boldsymbol{\xi}}{(1+\epsilon)\underline{x}}$$
$$= \underline{\omega} \left( K\theta + S_r p_w \right) \frac{3\underline{x}}{m_v} \frac{\underline{\boldsymbol{Y}}}{\underline{y}}. \quad (33)$$

By (33) we demonstrate that under small isotropic deformation the total force state obtained from the correspondence material model is equivalent to the one obtained from the ordinary material model. With (30) the nonlocal pressure gradient and water flow state can be also written as follows.

$$\widetilde{\nabla \Phi} = \frac{3}{m_v} \left( \int_{\mathcal{H}} \underline{\omega} \, \underline{\Phi} \boldsymbol{\xi} \, dV' \right) \mathbf{1} = \frac{3}{m_v} \int_{\mathcal{H}} \underline{\omega} \, \underline{\Phi} \boldsymbol{\xi} \, dV' \quad (34),$$

$$\underline{Q} = \frac{3}{m_v} \underline{\omega} \rho_w \boldsymbol{q} (\mathbf{1} \boldsymbol{\xi}) = \frac{3}{m_v} \underline{\omega} \rho_w \boldsymbol{q} \boldsymbol{\xi}. \quad (35)$$

*2.2. Fracture unsaturated periporomechanics*

In this section, we present fracture unsaturated periporomechanics that extends unsaturated periporomechanics to model fracture in unsaturated porous media. In fracture periporomechanics, both bulk and fracture space are represented by mixed material points that have two types of degree of freedom, i.e., displacement and fluid pressure. The fracturing process in periporomechanics is modeled following the bond-breakage concepts in the original peridynamics for solids. Figure 2 schematically represents the bond breakage concept for modeling fracture formation in unsaturated periporomechanics. The broken poromechanical bond will not be considered when determining the effective force state at material point *x* through the constitutive model. The effective force state is zero on the broken poromechanical bond while pore fluid pressure remains. To facilitate the modeling of unsaturated fluid flow in fracture space, the mixed material point in fracture space is named fracture point that has two fluid pressures, i.e., bulk fluid pressure $p_w$ and fracture fluid pressure $p_f$.

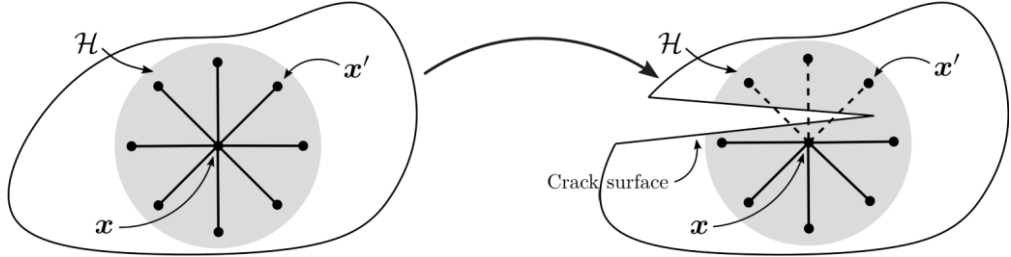

Figure 2: Conceptual sketch of fracturing in periporomechanics: (a) a material point in a porous body, and (b) a material point with bond breakage in a fractured porous body.

*2.2.1. Fracture criterion based on effective force state*

The bond breakage can be modeled either by the critical stretch criterion or energy criterion (e.g.,[18, 40]). In unsaturated porous media, the bond-breakage can be triggered by deformation or matric suction in unsaturated porous media as illustrated by the energy based criterion in what follows. In this study, the bond-breakage criteria in unsaturated periporomechanics are based on the deformation energy stored in a poromechanical bond. The effective force state that is an energy conjugate of the deformation state [12] is used to determine the deformation energy, which is different from the single phase peridynamics analysis. The effective force principle captures the effect of pressure developed in the pore fluid on the mechanical response of the skeleton. It is assumed that the energy density stored in a solid bond, $\varpi$, is fully recoverable until it exceeds some critical value $\varpi_{cr}$. With effective force state the energy density in a solid bond $\xi$ reads

$$\underline{\varpi} = \int_0^{t_l} \left(\underline{\overline{T}} - \underline{\overline{T}}'\right) \dot{\underline{\eta}} \, \mathrm{d}t = \int_0^{t_l} \left[(\underline{T} + S_r \underline{T}_w) - (\underline{T}' + S_r' \underline{T}_w')\right] \dot{\underline{\eta}} \, \mathrm{d}t, \tag{36}$$

where $\eta = u^0 - u$ is the relative displacement vector and $t_l$ is the total loading time. It is noted that the energy-based criterion of (36) incorporates the coupling effect of skeleton deformation, matric suction, and degree of saturation on the bond breakage. In fracture unsaturated periporomechanics, bond breakage is modeled through the influence function at the material constitutive level for the solid and fluid phases. The influence function will be replaced by a new influence function $\underline{\varrho \omega}$, where $\underline{\varrho}$ is defined as

$$\underline{\varrho} = \begin{cases} 0 & , \text{ if } \underline{\varpi} \geq \varpi_{cr} \\ 1 & , \text{ otherwise} \end{cases} \tag{37}$$

In periporomechanics, with the effective force concept the failure of solid skeleton is modeled through a scalar damage variable $\varphi$ that tracks the progressive failure of unsaturated porous media. This damage variable is defined as the fraction of broken solid bonds at a material point in its horizon

$$\varphi = 1 - \frac{\int_\mathcal{H} \underline{\varrho \omega} \, \mathrm{d}V'}{\omega_0}, \tag{38}$$

where $\varphi \in [0,1]$ and $\omega_0$ is defined by (8). It is hypothesized that once a solid bond breaks it will not sustain any mechanical load and the load at a material point will re-distribute to unbroken bonds. If enough bonds break (i.e., $\varphi >= \varphi_{cr}$) and coalesce into a surface, fracture will form and propagate naturally. In this study, we assume that cracks can be identified if $\varphi > 0.5$.

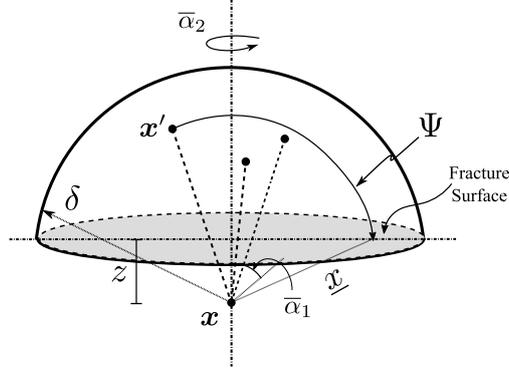

Figure 3: Schematic representation of the spherical coordinate system for the evaluation of the energy per unit fracture area $G_c$ for complete separation of two halves of a porous body (adapted from [24]).

The critical energy per unit fracture area, $G_c$, can be determined from $\varpi_{cr}$ through a spherical coordinate system as sketched in Figure 3. Referring to Figure 3, the energy per unit fracture area $G_c$ for completely separating the body into two halves reads

$$G_c = \int_0^\delta \int_0^{2\pi} \int_z^\delta \int_0^\Psi (\varpi_{cr})\underline{x}^2 \sin\overline{\alpha}_1 \, d\overline{\alpha}_1 \, d\underline{x} \, d\overline{\alpha}_2 \, dz,$$
$$= \frac{\pi\delta^4}{4}\varpi_{cr}, \tag{39}$$

where $\Psi = \cos^{-1}(z/\underline{x})$. It follows from (39) that the critical energy density for bond breakage can be written as

$$\varpi_{cr} = \frac{4G_c}{\pi\delta^4}. \tag{40}$$

From linear elastic fracture mechanics [23] for mode I fracture $G_c$ reads

$$G_c = K_I^2(1-\nu^2)/E, \tag{41}$$

where $E$ is the Young's modulus, $\nu$ is the Poisson's ratio and $K_I$ is the fracture toughness of mode I. Then from (41) the critical energy density for mode I fracture can be written as

$$\varpi_{cr} = \frac{4K_I^2}{\pi E \delta^4}(1-\nu^2). \tag{42}$$

It follows from (42) that the critical energy density for bond breakage can be a material property of solid phase. Following the above lines, the critical energy density can be determined for modes II and III fracture.

We assume that a fracture space is formed between two adjacent material points at $x$ and $x'$ if both $\varpi \geq \varpi_{cr}$ as well as $\varphi \geq \varphi_{cr}$ and $\varphi' \geq \varphi_{cr}$ at $x$ and $x'$, respectively. Both the material points are defined as the fracture mixed point that represent both the bulk and fractured space.

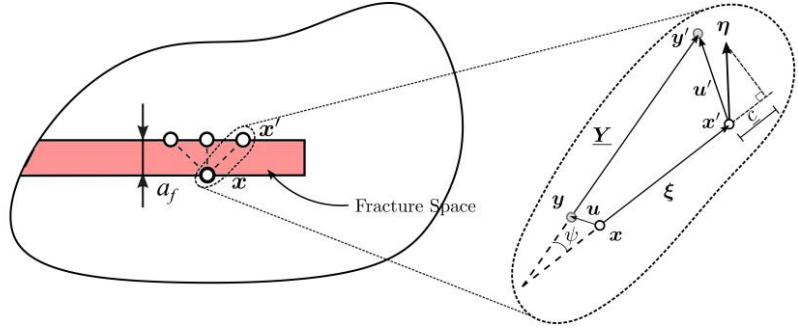

Figure 4: Schematic representation of the kinematics of a bond across fracture.

It is further assumed that the fracture mixed point has two degrees of freedom of fluid pressure, i.e., the bulk fluid pressure and fracture fluid pressure, which will be used to model unsaturated fluid flow in fractured space in Section 2.2.3. Figure 4 draws the kinematics of a bond across fracture. The relative deformation $\boldsymbol{\eta}$ can be decomposed into two components as sketched in Figure 4. The two components can be assumed to represent the opening displacement and the dislocation of the crack, respectively. To determine the crack width, we define $\underline{c}$ as the crack aperture related to the opening displacement as proposed in [34],

$$\underline{c} = \underline{y}\cos\psi - \underline{x}, \qquad (43)$$

where $\psi$ is defined in Figure 4. Therefore, the crack width at the fracture point *x* can be approximated by the average of bond apertures of all broken bonds

$$a_f = \frac{\int_{\mathcal{H}} \overline{\underline{\varrho}}\,\underline{c}\,\mathrm{d}V'}{\int_{\mathcal{H}} \overline{\underline{\varrho}}\,\mathrm{d}V'}, \qquad (44)$$

where $\overline{\underline{\varrho}} = 1- \underline{\varrho}$. Next, we present the equation of motion and mass balance of fracturing unsaturated porous media.

*2.2.2. Equation of motion*

In unsaturated periporomechanics, through the effective force concept the same equation of motion are applied on and off the crack surface or crack tip of the solid skeleton in unsaturated porous media. This is an advantage for modeling fracture in unsaturated porous media. For the fully coupled processes of solid deformation and unsaturated fluid flow the equation of motion of a fracture point is different from that of a bulk material point at the constitutive model level. For a bulk material point, the motion of equation can be written by (1). To incorporate bond-breakage at a bulk material point $x$, the effective force state and fluid force state along $\xi$ are written as respectively

$$\underline{\overline{T}} = \varrho\underline{\omega}\left(\overline{P}K^{-1}\xi + \frac{GC}{\omega_0}\underline{\mathscr{R}}^s\right), \qquad (45)$$

$$\underline{\overline{T}}_w = -\underline{\omega}P_w K^{-1}\xi, \qquad (46)$$

where,

$$P_w = Jp\mathbf{1}\widetilde{F}^{-T}. \qquad (47)$$

It follows from (45) and (46) that the total force state along a broken bond is

$$\underline{T} = -\underline{\omega}S_r P_w K^{-1}\xi. \qquad (48)$$

The broken bonds are eliminated through $\varrho$ at the constitutive model level (e.g., the nonlocal deformation gradient and the effective force state) [25, 41, 58]. As assumed, along a broken bond the effective force state vanishes. The total force state only consists of the fluid force state that can be written as

$$\underline{\mathbf{T}}_w = -\frac{3\omega}{m_v} p \frac{\underline{\mathbf{Y}}}{\underline{y}}, \tag{49}$$

where the ordinary fluid force state is used on broken bonds for simplicity.

For a fracture point, similar to (1) the equation of motion can be written as

$$\rho \ddot{\mathbf{u}} = \int_{\mathcal{H}} \left( \underline{\mathbf{T}} - \underline{\mathbf{T}}' \right) \, \mathrm{d}V' - \int_{\mathcal{H}} \left( S_l \underline{\mathbf{T}}_l - S_l' \underline{\mathbf{T}}_l' \right) \, \mathrm{d}V' + \rho \mathbf{g}, \tag{50}$$

where

$$S_l \underline{\mathbf{T}}_l = \begin{cases} S_{r,f} \underline{\mathbf{T}}_f & \text{if } \varphi \ \& \ \varphi' \geq \varphi_{\mathrm{cr}} \\ S_r \underline{\mathbf{T}}_w & \text{otherwises} \end{cases}. \tag{51}$$

In (51), $\underline{\mathbf{T}}_f$ is the fluid force state between two fracture points at $x$ and $x'$

$$\underline{\mathbf{T}}_f = -\frac{3\omega}{m_v} p_f \frac{\underline{\mathbf{Y}}}{\underline{y}}. \tag{52}$$

where $p_f$ and $S_{r,f}$ are fracture fluid pressure and degree of saturation respectively at fracture point $x$.

### 2.2.3. Balance of mass

The balance of mass in the continuous porous space (Figure 6(a)) can be written by (3) at both bulk and fracture material points. In (3), the broken bond (i.e., $\varrho$ = 0) associated with $x$ will not be considered in the second and third terms and at the constitutive model level by multiplying $\underline{\omega}$ by %. For instance, equations (34) and (35) can be rewritten as

$$\widetilde{\nabla \Phi} = \frac{3}{m_v} \int_{\mathcal{H}} \varrho \omega \, \Phi \boldsymbol{\xi} \, \mathrm{d}V' \tag{53},$$

$$\underline{Q} = \frac{3}{m_v} \varrho \omega \rho_w \mathbf{q} \boldsymbol{\xi}. \tag{54}$$

The sink term will be null for a bulk material point. For a fracture material point, the sink term in (3) can be determined assuming that the fluid flow from the pore space into the fracture space follows the generalized Darcy's law for unsaturated fluid flow along the direction normal to the fracture surface [6, 31]. It follows from this assumption that sink/source term $Q_s$ for the fluid flow from the bulk into the fracture at $x$ can be written as

$$Q_s = A \left[ -\frac{k^r k_w}{\mu_w} \left( \frac{p_f - p_w}{l_x} \right) \right] / V, \tag{55}$$

where $p_w$ and $p_f$ are water pressures in the bulk and fracture space respectively, $A$ and $V$ are the cross-sectional area and volume of a material point assuming uniform spatial discretization, $l_x$ = $d/2$ and $d$ is the edge dimension of a cubic material point. Given the source term, following the formulation for unsaturated fluid flow in porous space the mass balance equation of unsaturated fluid flow in fracture space can be written as

$$\frac{\partial S_{r,f}}{\partial t} + \frac{1}{\rho_w} \int_{\mathcal{H}} \left( \underline{Q}_f - \underline{Q}'_f \right) \mathrm{d}V' - Q_s = 0, \tag{56}$$

where $\underline{Q}_f$ and $\underline{Q}'_f$ are the fluid flow states at fracture points at $x$ and $x'$, respectively, $S_{r,f}$ is the degree of saturation in fracture space that can be determined by the soil-water retention curve in (16). In (56), it is assumed that in fracture space $\phi = 1$ and the volume coupling term vanishes [6]. The fracture flow state $Q_f$ is determined by (35) as

$$\underline{Q}_f = \frac{3}{m_v} \omega \rho_w \boldsymbol{q}_f \boldsymbol{\xi} \tag{57}$$

where $q_f$ is the fluid flow vector in fracture space. Through Darcy's law for unsaturated fluid flow [3, 13] the fracture fluid flow vector $q_f$ can be written as

$$\boldsymbol{q}_f = -\frac{k_f^r k_f}{\mu_w} \widetilde{\boldsymbol{\nabla} \Phi}_f, \tag{58}$$

where $k_f^r$ is the relative permeability, $k_f$ is the intrinsic permeability of fracture space, and $\widetilde{\nabla \Phi}_f$ is the nonlocal fracture fluid pressure gradient.

$$k_f^r = S_{r,f}^{1/2} \left[ 1 - (1 - S_{r,f}^{1/m})^m \right]^2. \tag{59}$$

Given the fracture width $a_f$ in (44), $k_f$ can be estimated through the celebrated cubic law (e.g.,[6])

$$k_f = \frac{a_f^2}{12}. \tag{60}$$

The nonlocal gradient of fracture fluid pressure can be determined

$$\widetilde{\nabla \Phi}_f = \frac{3}{m_v} \int_{\mathcal{H}} \underline{\omega} \, \underline{\Phi}_f \boldsymbol{\xi} \, \mathrm{d}V', \tag{61}$$

where

$$\underline{\Phi}_f = p'_f - p_f, \tag{62}$$

and $p_f$ and $p'_f$ are fracture fluid pressures at fracture points $x$ and $x'$, respectively.

In summary in the proposed fracturing unsaturated periporomechanics the fundamental unknowns are displacement and fluid pressures in the bulk and fracture space. For a bulk material point, the governing equations consist of the equation of motion and the mass balance equation, i.e., (1) and (3). For a fracture material point, the governing equations are the equation of motion, the mass balance equation in the bulk and the mass balance equation of fluid flow in fracture space, i.e., (50), (3), and (56). In the following section, we present an implicit fractional step algorithm to solve the coupled governing equations of fracturing unsaturated periporomechanics.

## 3. Numerical implementation

### 3.1. Fractional step method

The formulated unsaturated fracturing periporomechanics is numerically implemented through the fractional step/staggered algorithm with a two-state undrained operator split [15, 16, 47, 49] and mixed meshfree method [14]. The two-stage operator split converts the coupled fracturing periporomechanics problem into an undrained mechanical problem and an unsaturated fluid flow problem in the deformed solid configuration. We note the drained operator split that involves simply freezing the fluid pressure during the mechanical stage is only conditionally stable [48, 59]. The undrained operator split preserves the contractivity of the original coupled problem and is unconditionally stable [49]. It is noted that both monolithic and staggered approach can be exploited to solve the fracturing unsaturated periporomechanics model in this study. In the monolithic approach the full system of coupled algebraic equations are solved simultaneously at the same time step. With an implicit time integration scheme, it is unconditionally stable and preserves the strong physical coupling between the phases. However, it can lead to large complex systems that may be asymmetric and cannot leverage the different phenomenological time scales of individual phases [15]. We refer to the distinguished literature (e.g., [15, 16]) for an in-depth discussion of the monolithic and fractional step/staggered algorithms for numerically implementing coupled multiphysics frameworks for porous media.

Figure 5 provides a flowchart of the proposed implicit fractional step formulation and Algorithm 1 presents a detailed step-by-step procedure.

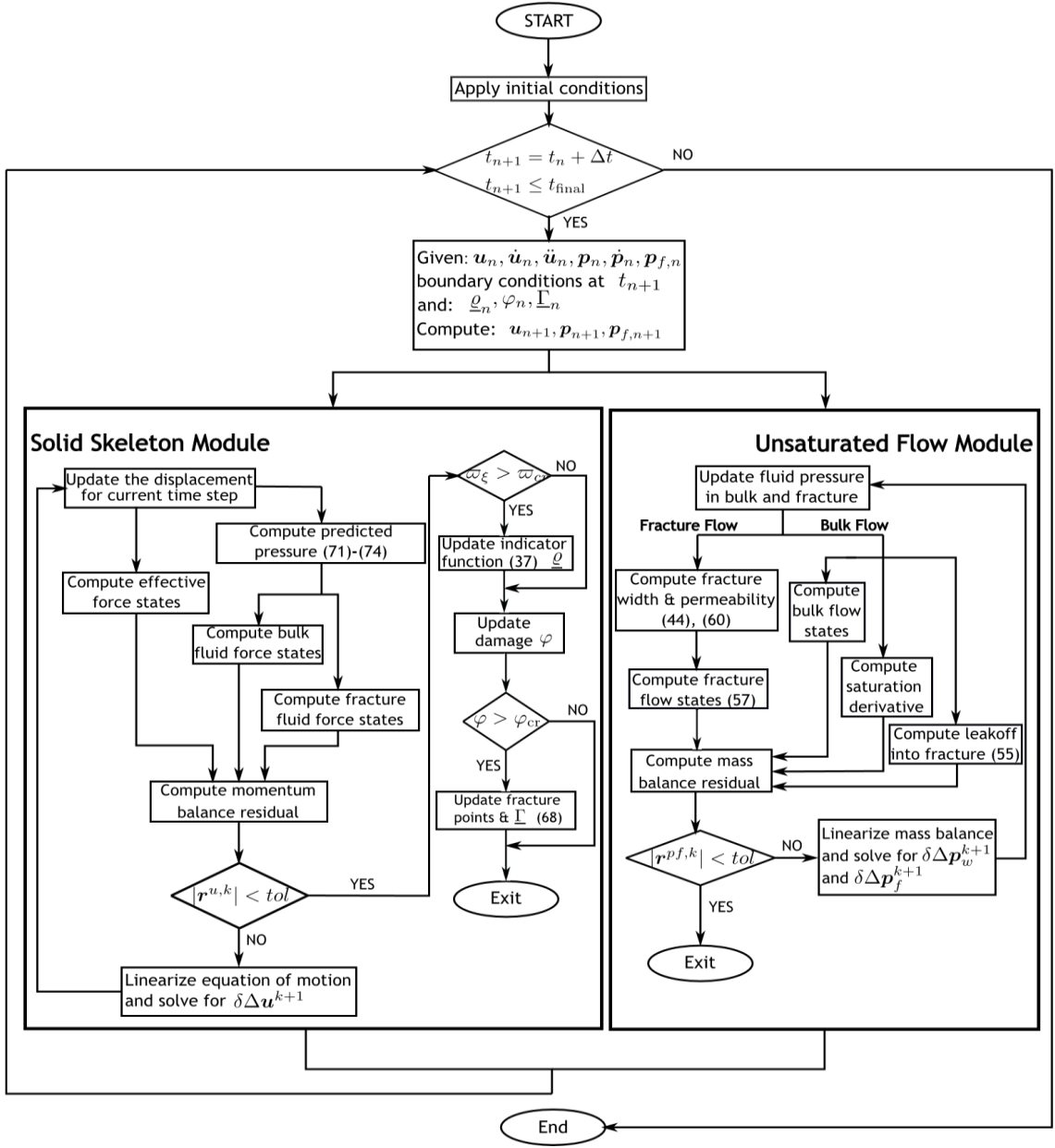

Figure 5: Flowchart for the implicit fractional step algorithm for fracturing unsaturated periporomechanics.

*3.1.1. Time discretization of the equation of motion under undrained condition*

In the mechanical stage, the motion of equation will be solved first under undrained condition through Newton's method. At time step $n$, $u_n$, $\dot{u}_n$, $\ddot{u}_n$ are known. Let $\Delta u_{n+1} = u_{n+1} - u_n$, in line with Newmark's method [60] the displacement, velocity, and accelerations of the solid skeleton can be written as

$$u_{n+1} = u_n + \Delta u_{n+1}, \tag{63}$$

$$\dot{u}_{n+1} = \frac{2\beta_2}{\beta_1 \Delta t}\Delta u_{n+1} - \left(\frac{2\beta_2}{\beta_1} - 1\right)\dot{u}_n - \Delta t\left(\frac{\beta_2}{\beta_1} - 1\right)\ddot{u}_n, \tag{64}$$

$$\ddot{u}_{n+1} = \frac{2}{\beta_1 \Delta t^2}\Delta u_{n+1} - \frac{2}{\beta_1 \Delta t}\dot{u}_n - \left(\frac{1}{\beta_1} - 1\right)\ddot{u}_n, \tag{65}$$

where $\beta_1$ and $\beta_2$ are numerical integration parameters. For unconditional stability [16] $\beta_1 \geq \beta_2 \geq 0.5$.

---

**Algorithm 1** Global fractional-step time integration scheme for fractured porous media

1: **procedure** GIVEN $u_n, p_n, p_{f,n}, t_n$ AND $\Delta t$, COMPUTE $u, p$ AND $p_f$
2:     $t = t_n + \Delta t$
3:     **while** $t \leq t_{\text{final}}$ **do**
4:         Apply boundary conditions
5:         Compute the effective force and fluid force via Algorithm 2.
6:         Compute balance of momentum residual $\mathcal{R}^{u,0}$
7:         Set $k = 0$ and tol $= 10^{-6}$
8:         **if** $\mathcal{R}^{u,0} >$ tol **then**
9:             **while** $|\mathcal{R}^{u,k}|/|\mathcal{R}^{u,0}| >$ tol **do**
10:                 Construct tangent operator for balance of momentum $\mathcal{A}^{u,k} = \left[\partial \mathcal{R}^{u,k}/\partial \Delta u^k\right]$
11:                 Solve the linear system $\mathcal{A}^{u,k}\delta\Delta u^{k+1} = -\mathcal{R}^{u,k}$ for $\delta\Delta u^{k+1}$
12:                 Update $u^{k+1}, \dot{u}^{k+1}$ and $\ddot{u}^{k+1}$ using (63), (64), (65)
13:                 Update the effective force and fluid force via Algorithm 2
14:                 Update the residual $\mathcal{R}^{u,k+1}$
15:                 Set $k \leftarrow k+1$
16:             **end while**
17:         **end if**
18:         Update list of broken bonds and fracture points via Algorithm 3
19:         Compute the bulk flow and fracture flow via Algorithm 4
20:         Compute the balance of mass residual $\mathcal{R}^{pf,0} = \{\mathcal{R}^{p,0}, \mathcal{R}^{f,0}\}$
21:         Set $k = 0$
22:         **if** $\mathcal{R}^{pf,0} >$ tol **then**
23:             **while** $|\mathcal{R}^{pf,k}|/|\mathcal{R}^{pf,0}| >$ tol **do**
24:                 Construct tangent $\mathcal{A}^{pf,k} = \left[\{\partial \mathcal{R}^{pf,k}/\partial \Delta p^k\}, \{\partial \mathcal{R}^{pf,k}/\partial \Delta p_f^k\}\right]$
25:                 Solve the linear system $\mathcal{A}^{pf,k}\delta\Delta p^{k+1} = -\{\mathcal{R}^{p,k}, \mathcal{R}^{f,k}\}$ for $\delta\Delta p^{k+1}$
26:                 Update $p^{k+1}$ and $p_f^{k+1}$ using (80) and (81)
27:                 Update the fracture and bulk flow via Algorithm 4
28:                 Update the residual $\mathcal{R}^{pf,k+1}$
29:                 Set $k \leftarrow k+1$
30:             **end while**
31:         **end if**
32:         Update $u_n \leftarrow u, p_n \leftarrow p, p_{f,n} \leftarrow p_f$
33:     **end while**
34: **end procedure**

---

Figure 6 presents a schematic of (a) a material point in the bulk and (b) a material point adjacent to the fracture. At $t_{n+1}$, the residual of the equation of motion at a bulk point reads,

$$r_{n+1}^u = \rho\ddot{u} - \int_\mathcal{H} \underline{\varrho}_n \left(\underline{\boldsymbol{T}} - \underline{\boldsymbol{T}}'\right) dV' + \int_\mathcal{H} \left(S_r \underline{\boldsymbol{T}}_w - S_r' \underline{\boldsymbol{T}}_w'\right) dV' - \rho\boldsymbol{g}. \tag{66}$$

For a fracture point, the residual for equation of motion (from (50)) is expressed as

$$r^u_{n+1} = \rho\ddot{\boldsymbol{u}} - \int_{\mathcal{H}} \underline{\varrho}_n \left(\overline{\boldsymbol{T}} - \overline{\boldsymbol{T}}'\right)\, \mathrm{d}V' + \int_{\mathcal{H}} (1-\underline{\Gamma}_n)\left(S_w\underline{\boldsymbol{T}}_w - S'_w\underline{\boldsymbol{T}}'_w\right)^k\, \mathrm{d}V'$$
$$+ \int_{\mathcal{H}} \underline{\Gamma}_n \left(S_f\underline{\boldsymbol{T}}_f - S'_f\underline{\boldsymbol{T}}'_f\right)\, \mathrm{d}V' - \rho^k\boldsymbol{g}, \qquad (67)$$

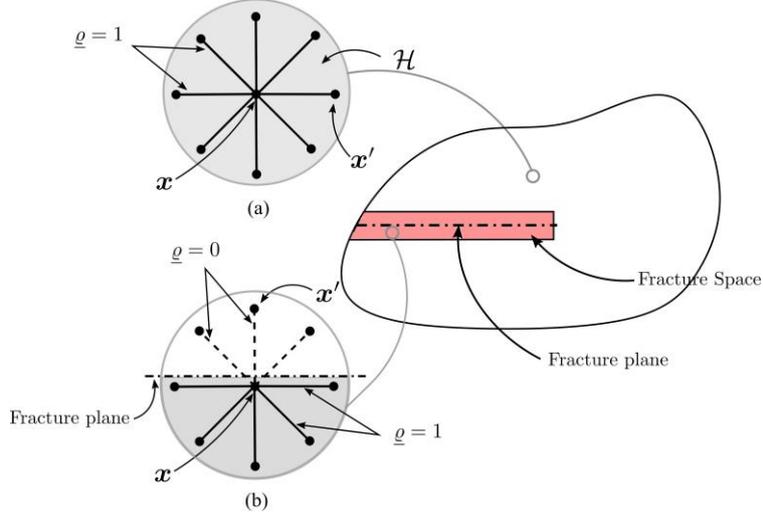

Figure 6: Schematic representation of (a) a point in the bulk (b) a fracture point adjacent to the crack.

where

$$\underline{\Gamma}_n = \begin{cases} 1, & \text{if } \boldsymbol{x} \ \& \ \boldsymbol{x}' \text{ are fracture points at } t_n, \\ 0, & \text{otherwise.} \end{cases} \qquad (68)$$

---
**Algorithm 2** Compute the effective force and fluid force
---
1: **procedure** GIVEN $\boldsymbol{u}^k, \dot{\boldsymbol{u}}^k, \ddot{\boldsymbol{u}}^k, \underline{\varrho}_n, \underline{\Gamma}_n$ CONSTRUCT EFFECTIVE FORCE VECTOR $\mathscr{T}^k$ AND FLUID FORCE VECTOR $\mathscr{T}^k_w$
2:     **for** all points **do**
3:         Compute weighted volume $m_v$ using (20) or shape tensor $\boldsymbol{K}$ using (29)
4:         Compute deformation gradient $\widetilde{\boldsymbol{F}}^k$ using (10)
5:         Compute effective Piola-Kirchhoff stress $\overline{\boldsymbol{P}}^k$ using constitutive model
6:         Compute volume coupling term $\mathbb{V}^k$ using (99)
7:         Compute force states $\overline{\underline{\boldsymbol{T}}}^k$, $\underline{\boldsymbol{T}}^k_w$ and $\underline{\boldsymbol{T}}^k_f$ using (45), (46), and (52)
8:         Compute the rate fluid pressure $\dot{\widetilde{p}}^k_w$ using (71)
9:         Compute the rate fracture fluid pressure $\dot{\widetilde{p}}^k_f$ using (72)
10:        Compute bulk pressure $\widetilde{p}_w$ using (73) and fracture pressure $\widetilde{p}_f$ using (74)
11:        Compute effective force $\mathscr{T}^k$ using (93)
12:        Compute fluid force $\mathscr{T}^k_w$ using (94)
13:     **end for**
14: **end procedure**

The fluid pressure terms in (66) and (67) are determined from the explicit fluid pressure predictors under undrained condition that are defined in what follows. At time step $n$, $p_{w,n}$ and $p_{f,n}$, and their rate forms, $\dot{p}_{w,n}$ and $\dot{p}_{f,n}$ are known. Let $\tilde{p}_{w,n+1}$ and $\tilde{p}_{f,n+1}$ be the explicit fluid pressure predictors at time step $n+1$ under undrained condition. It follows from Newmark's method [14, 16] that the rate forms $\dot{\tilde{p}}_{w,n+1}$ and $\dot{\tilde{p}}_{f,n+1}$ can be written as

$$\dot{\tilde{p}}_{w,n+1} = \frac{1}{\beta_3 \Delta t} \Delta \tilde{p}_{w,n+1} - \left(\frac{1}{\beta_3} - 1\right) \dot{p}_{w,n} \tag{69}$$

$$\dot{\tilde{p}}_{f,n+1} = \frac{1}{\beta_3 \Delta t} \Delta \tilde{p}_{f,n+1} - \left(\frac{1}{\beta_3} - 1\right) \dot{p}_{f,n}, \tag{70}$$

where $\Delta \tilde{p}_{w,n+1} = \tilde{p}_{w,n+1} - p_{w,n}$ and $\Delta \tilde{p}_{f,n+1} = \tilde{p}_{f,n+1} - p_{f,n}$, and $\beta_3$ is an integration parameter ($\beta_3 \geq 0.5$ for unconditional stability). Under undrained condition, $\tilde{p}_{w,n+1}$ and $\tilde{p}_{f,n+1}$ can be determined from the mass balance equations of unsaturated fluid flow in the bulk and fracture as follows.

$$\dot{\tilde{p}}_{w,n+1} = -\left(\phi \frac{\partial S_{r,n}}{\partial p_{w,n}}\right)^{-1} \left[ S_{r,n} \int_{\mathcal{H}} \underline{\varrho}_n \left(\underline{\dot{Y}} - \underline{\dot{Y}}'\right) dV' + \frac{1}{\rho_w} \int_{\mathcal{H}} \underline{\varrho}_n \left(\underline{Q} - \underline{Q}'\right)|_n dV' + \underline{Q}_{s,n} \right] \tag{71}$$

$$\dot{\tilde{p}}_{f,n+1} = -\left(\frac{\partial S_{r,f}}{\partial p_f}\right)_n^{-1} \left[ \frac{1}{\rho_w} \int_{\mathcal{H}} \underline{\Gamma}_n \left(\underline{Q}_f - \underline{Q}'_f\right)_n dV' - \underline{Q}_{s,n} \right]$$
$$= \dot{p}_{f,n}. \tag{72}$$

We note that $\underline{Q}_{s,n} = 0$ in (71) at a mixed point in the bulk. Thus, (71) represents the rates of water pressure predictors at the bulk and fracture points. Then, it follows from (69), (70), (71), and (72) that the explicit fluid pressure predictors can be written as

$$\tilde{p}_{w,n+1} = \beta_3 \Delta t \dot{\tilde{p}}_{w,n+1} + (1-\beta_3) \Delta t \dot{p}_{w,n} + p_{w,n}, \tag{73}$$

$$\tilde{p}_{f,n+1} = \beta_3 \Delta t \dot{\tilde{p}}_{f,n+1} + (1-\beta_3) \Delta t \dot{p}_{f,n} + p_{f,n}$$
$$= \Delta t \dot{p}_{f,n} + p_{f,n}. \tag{74}$$

The fluid pressure predictors from (73) and (74) are then used to determine the fluid force states in (66) and (67). Substituting (63), (64), and (65) into (66) and (67) and using (73) and (74), $\Delta u_{n+1}$ can be solved by Newton's method in what follows (see Algorithm 2). Let $k$ be the iteration number

$$\boldsymbol{\mathcal{R}}^{u,k+1} = \boldsymbol{\mathcal{R}}^{u,k} + \boldsymbol{\mathcal{A}}^{u,k} \delta \Delta \boldsymbol{u}^{k+1} \approx 0, \tag{75}$$

where $R^u$ is the global residual vector of displacements and $A^u$ is the global tangent operator of the motion equation

$$\boldsymbol{\mathcal{A}}^{u,k} = \frac{\partial \boldsymbol{\mathcal{R}}^{u,k}}{\partial \Delta \boldsymbol{u}^k}. \tag{76}$$

Solving (75) we obtain,

$$\delta \Delta \boldsymbol{u}^{k+1} = -(\boldsymbol{\mathcal{A}}^k)^{-1} \boldsymbol{\mathcal{R}}^{u,k}. \tag{77}$$

Finally, we have,

$$\Delta \boldsymbol{u}^{k+1} = \Delta \boldsymbol{u}^k + \delta \Delta \boldsymbol{u}^{k+1}. \tag{78}$$

At the end of mechanical stage, the broken bonds and indicator functions (i.e., (68)) are updated following the Algorithm 3. The energy dissipation of individual bonds is written as

$$\underline{\varpi}_{n+1} = \underline{\varpi}_n + \left(\underline{T} - \underline{T}'\right)_{n+1} \Delta \eta$$
$$= \underline{\varpi}_n + \left(\underline{T} - \underline{T}'\right)_{n+1} (\Delta u' - \Delta u). \quad (79)$$

The database of damaged bonds is updated for bonds with $\varpi > \varpi_{cr}$ via (37) and then $\varphi$ is updated at each material point. Using $\varphi$, the database of fracture points is updated.

The solution of the mechanical stage and the rate forms of the fluid pressure predictors at time $t_{n+1}$ from the mechanical solver are passed to the unsaturated fluid flow solver for the fluid flow stage at time step $t_{n+1}$.

---

**Algorithm 3** Update list of broken bonds and fracture points

1: **procedure** GIVEN $\underline{T}$ AND $\Delta \eta$ UPDATE INFLUENCE FUNCTION $\varrho$ AND DAMAGE $\varphi$
2:   **for** all points **do**
3:     **for** each neighbor **do**
4:       Compute bond energy $\underline{\varpi}$ using (79)
5:       **if** $\underline{\varpi} > \varpi_{cr}$ **then**
6:         Update $\varrho$ using (37) to reflect bond damage status
7:         Update damage $\varphi$ using (38)
8:         Sum the energy in the bond into the total energy dissipated at the point
9:       **end if**
10:    **end for**
11:    **if** $\varphi \geq \varphi_{cr}$ **then**
12:      set $x$ = fracture point
13:    **end if**
14:  **end for**
15: **end procedure**

---

*3.1.2. Time discretization of the mass balance equation in the deformed configuration*

Given $p_{w,n}$, $p_{f,n}$, $\dot{p}_{w,n}$, $\dot{p}_{f,n}$ and $\dot{u}_{n+1}$, the unsaturated fluid flow stage solves $p_{w,n+1}$ and $p_{f,n+1}$ in the updated solid configuration (i.e., $u_{n+1}$) at $t_{n+1}$ using Newton's method. Let $\Delta p_{w,n+1} = p_{w,n+1} - p_{w,n}$ and $\Delta p_{f,n+1} = p_{f,n+1} - p_{f,n}$. It follows from Newmark's method [60] that

$$p_{w,n+1} = p_{w,n} + \Delta p_{w,n+1}, \quad (80)$$
$$p_{f,n+1} = p_{f,n} + \Delta p_{f,n+1}, \quad (81)$$
$$\dot{p}_{w,n+1} = \frac{1}{\beta_3 \Delta t} \Delta p_{w,n+1} - \left(\frac{1}{\beta_3} - 1\right) \dot{p}_{w,n}, \quad (82)$$
$$\dot{p}_{f,n+1} = \frac{1}{\beta_3 \Delta t} \Delta p_{f,n+1} - \left(\frac{1}{\beta_3} - 1\right) \dot{p}_{f,n}, \quad (83)$$

At $t_{n+1}$, the residuals of the mass balance equations can be written as

$$r^p_{n+1} = \phi \left( \frac{\partial S_r}{\partial p_w} \right) \dot{p}_w + \frac{1}{\rho_w} \int_{\mathcal{H}} \varrho \left( \underline{Q} - \underline{Q}' \right) \, \mathrm{d}V' + S^k_r \int_{\mathcal{H}} \varrho \left( \underline{\dot{y}} - \underline{\dot{y}}' \right) \, \mathrm{d}V' + \mathcal{Q}_s \tag{84}$$

$$r^f_{n+1} = \left( \frac{\partial S_{r,f}}{\partial p_f} \right) \dot{p}_f + \frac{1}{\rho_w} \int_{\mathcal{H}} \Gamma \left( \underline{Q}_f - \underline{Q}'_f \right) \mathrm{d}V' - \mathcal{Q}_s. \tag{85}$$

Substituting (80), (81), (82) and (83) into (84) and (85), $\Delta p_{w,n+1}$ and $\Delta p_{f,n+1}$ can be solved by Newton's method in what follows (see Algorithm 4 for details). Again, let $k$ be the iteration number

$$\begin{Bmatrix} \mathcal{R}^{p,k+1} \\ \mathcal{R}^{f,k+1} \end{Bmatrix} = \begin{Bmatrix} \mathcal{R}^{p,k} \\ \mathcal{R}^{f,k} \end{Bmatrix} + \mathcal{A}^{pf,k} \begin{Bmatrix} \delta \Delta p_w^{k+1} \\ \delta \Delta p_f^{k+1} \end{Bmatrix} \approx 0, \tag{86}$$

where $\mathcal{R}^p$ and $\mathcal{R}^f$ are the global residual vectors of the mass balance equations of the bulk and fracture respectively, $\Delta p_w$ and $\Delta p_f$ are the global fluid pressure increments in the bulk and fracture respectively, and $\mathcal{A}^{pf}$ is the coupled tangent operator of the mass balance equations

$$\mathcal{A}^{pf,k} = \begin{Bmatrix} \frac{\partial \mathcal{R}^{p,k}}{\partial \Delta p_w^k} & \frac{\partial \mathcal{R}^{p,k}}{\partial \Delta p_f^k} \\ \frac{\partial \mathcal{R}^{f,k}}{\partial \Delta p_w^k} & \frac{\partial \mathcal{R}^{f,k}}{\partial \Delta p_f^k} \end{Bmatrix} \tag{87}$$

---

**Algorithm 4** Compute the unsaturated fluid flow in the bulk and fracture

---

1: **procedure** GIVEN $p^k, p_f^k, \dot{u}, \varrho, \Gamma$ CONSTRUCT BULK FLUID FLOW $\mathcal{Q}^k$ AND FRACTURE FLOW $\mathcal{Q}_f^k$
2:     **for** all points **do**
3:         Compute all pressure potential states $\underline{\Phi}^k$
4:         Compute weighted volume $m_v$ using (20) or the shape tensor $K$ using using (29)
5:         Compute pressure gradient $\widetilde{\nabla \Phi}^k$ using (14)
6:         Compute relative permeability $k^{r,k}$ using (17)
7:         Compute the flux vector $q^k$ using (14)
8:         Compute all flow states $\underline{Q}^k$ using (13)
9:         Compute fluid flow $\mathcal{Q}^k$ using (98)
10:     **end for**
11:     **for** all fracture points **do**
12:         Compute fracture width $a_f$ using (44)
13:         Compute fracture permeability $K_f$ using (60)
14:         Compute $\mathcal{Q}_s^k$ using (55)
15:         **for** each neighboring fracture point **do**
16:             Compute fracture pressure states $\underline{\Phi}_f^k$ using (62)
17:             Compute the fracture pressure gradient $\widetilde{\nabla \Phi}_f^k$ and fracture flux $q_f^k$ using (61) and (58)
18:             Calculate fracture flow states $\underline{Q}_f^k$ using (57)
19:         **end for**
20:         Compute the fracture flow $\mathcal{Q}_f^k$ using (102) and source term $\mathcal{Q}_s^k$
21:     **end for**
22: **end procedure**

---

Solving (86) gives

$$\begin{Bmatrix} \delta \Delta p_w^{k+1} \\ \delta \Delta p_f^{k+1} \end{Bmatrix} = -(\mathcal{A}^{pf,k})^{-1} \begin{Bmatrix} \mathcal{R}^{p,k} \\ \mathcal{R}^{f,k} \end{Bmatrix} \tag{88}$$

It is noted that $\mathring{p}_{w,n+1}$ and $\mathring{p}_{f,n+1}$ from the undrained mechanical stage are used as initial values of $p_{w,n+1}$ and $p_{f,n+1}$ (i.e., $k = 0$). Finally, we have,

$$\Delta p_w^{k+1} = \Delta p_w^k + \delta \Delta p_w^{k+1}, \tag{89}$$

$$\Delta p_f^{k+1} = \Delta p_f^k + \delta \Delta p_f^{k+1}. \tag{90}$$

*3.2. Mixed meshfree spatial discretization*

The formulated fracture periporomechanics model is discretized in space by a mixed Lagrangian-Eulerian meshfree scheme [14]. The uniform mixed material points are used in this study. Let $N_i$ be the number of mixed material points in the horizon of mixed material point $i$, P be the number of total mixed material points in the problem domain, and $A$ be the assembly operator [14].

*3.2.1. Spatial discretizion of the equation of motion*

It follows from (66) and (67) that the global residual vector of the equation of motion $\mathcal{R}_{n+1}^u$ at time step $n + 1$ can be constructed as [14],

$$\mathcal{R}_{n+1}^u = \mathcal{A}_{i=1}^{\mathcal{P}} \left( \mathcal{M}_i + \mathcal{T}_{s,i} - \mathcal{T}_{w,i} \right), \tag{91}$$

where $\mathcal{M}_i$, $\mathcal{T}_{s,i}$, and $\mathcal{T}_{w,i}$ are the inertial and gravity load, the effective force and the fluid force at mixed material point $i$.

$$\mathcal{M}_i = [\rho_s(1 - \phi_i) + \rho_w S_{r,i} \phi_i] (\ddot{u}_i - g) V_i, \tag{92}$$

$$\mathcal{T}_{s,i} = \sum_{j=1}^{\mathcal{N}_i} \varrho_{ij,n} \left( \overline{\boldsymbol{T}}_{ij} - \overline{\boldsymbol{T}}_{ji} \right) V_j V_i, \tag{93}$$

$$\mathcal{T}_{w,i} = \sum_{j=1}^{\mathcal{N}_i} (1 - \Gamma_{ij,n}) \left( S_{r,i} \underline{\boldsymbol{T}}_{w,ij} - S_{r,j} \underline{\boldsymbol{T}}_{w,ji} \right) V_j V_i$$

$$+ \sum_{j=1}^{\mathcal{N}_i} \Gamma_{ij,n} \left( S_{rf,i} \underline{\boldsymbol{T}}_{f,ij} - S_{rf,j} \underline{\boldsymbol{T}}_{f,ji} \right) V_j V_i, \tag{94}$$

where $j$ denotes material points in the horizon of $i$.

*3.2.2. Spatial discretization of the mass balance equation*

From (84) and (85), the global residual vectors of the mass balance equations $\mathcal{R}^p$ and $\mathcal{R}^f$ at time step $n + 1$ can be respectively written as

$$\mathcal{R}^p = \mathcal{A}_{i=1}^{\mathcal{P}} \left( \mathcal{X}_i + \mathcal{Q}_i + \mathcal{V}_i + \mathcal{Q}_{s,i} \right), \tag{95}$$

$$\mathcal{R}^f = \mathcal{A}_{i=1}^{\mathcal{P}_f} \left( \mathcal{X}_{f,i} + \mathcal{Q}_{f,i} - \mathcal{Q}_{s,i} \right), \tag{96}$$

where $P_f$ is the number of fracture points and

$$\mathcal{X}_i = -\Delta t \phi_i \frac{\partial S_{r,i}}{\partial s_i} \dot{p}_{w,i} V_i, \tag{97}$$

$$\mathcal{Q}_i = \frac{\Delta t}{\rho_w} \sum_{j=1}^{\mathcal{N}_i} \varrho_{ij,n+1} \left( \underline{Q}_{ij} - \underline{Q}_{ji} \right) V_j V_i, \tag{98}$$

$$\mathscr{V}_i = \Delta t S_{r,i} \sum_{j=1}^{\mathcal{N}_i} \underline{\varrho}_{ij,n+1} \left( \dot{\underline{\mathcal{V}}}_{ij} - \dot{\underline{\mathcal{V}}}_{ji} \right) V_j V_i, \tag{99}$$

$$\mathscr{Q}_{s,i} = \Delta t \left[ -\frac{k^r k_w}{\mu_w} \left( \frac{p_{f,i} - p_{w,i}}{l_i} \right) \right] A_i, \tag{100}$$

$$\mathscr{X}_{f,i} = -\Delta t \frac{\partial S_{rf,i}}{\partial p_f, i} \dot{p}_{f,i} V_i, \tag{101}$$

$$\mathscr{Q}_{f,i} = \frac{\Delta t}{\rho_w} \sum_{j=1}^{\mathcal{N}_i} \underline{\Gamma}_{ij,n+1} \left( \underline{Q}_{f,ij} - \underline{Q}_{f,ji} \right) V_j V_i. \tag{102}$$

## 4. Numerical examples

### 4.1. Example 1: Mode I crack propagation

In this example, we simulate mode I crack propagation and unsaturated fluid flow in a porous body under two-dimensions. Figure 7 depicts the problem geometry, boundary conditions and loading protocol. The pre-existing crack is modeled by eliminating interaction between material points across the crack plane [46]. All fluid phase boundaries of the specimen are assumed impermeable. The problem domain is discretized into 20000 mixed material points with $d$ = 2.5 mm. The solid phase is modeled using an isotropic elastic correspondence constitutive model [12, 13]. This applies to all subsequent examples in this section. It is assumed that $\rho_s$ = 2000 kg/m³, $\rho_w$ = 1000 kg/m³, $\mu_w$ = 1 × 10⁻³ Pa·s, initial porosity $\phi_0$ = 0.25 and $G_c$ = 95 J/m².

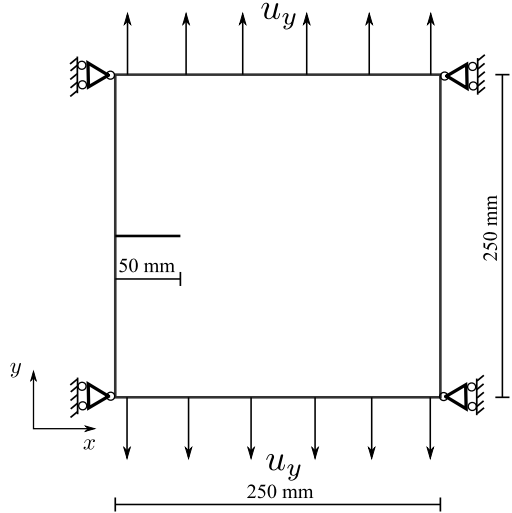

Figure 7: Problem setup for example 1.

Unless otherwise noted, these parameters are also used in all subsequent examples. The remaining material parameters in [61] are $K$ = 13.46 × 10⁶ kPa, $\mu_s$ = 10.95 × 10⁶ kPa, $k_w$ = 6 × 10⁻²¹ m², $n$ = 1.7844, $s_a$ = 18.6 × 10³ kPa. The horizon $\delta$ = 3.05$d$. The stabilization parameter $G$ = 1.0 [14]. As in [61], the specimen is prescribed with zero effective stress and pore water pressure. The loading rate is $\dot{u}_y$ = 2.35 ×10⁻⁵ mm/s. The total loading time $t$ = 2000 s with the time increment $\Delta t$ = 1 s. We investigate the influence of $\delta$ and $m = \delta/d$ on the crack propagation and unsaturated fluid flow. The loading curves are compared with the numerical results from the extended finite element method (XFEM) in [61].

### 4.1.1. Influence of $\delta$

The horizon $\delta$ plays a critical role in modeling the actual behavior of materials using peridynamics. However, its calibration for specific applications is still an open question in peridynamics. To examine the sensitivity of the results to the horizon $\delta$, we run the simulations with three values of horizon, $\delta$ = 11.3 mm, 7.5 mm, and 6 mm assuming $m$ = 3. The corresponding spatial discretizations consist of 9,000, 20,000, and 36,000 mixed points, respectively. All other parameters are identical. The results are presented in Figures 8, 9, and 10.

Figure 8 plots the loading curves from the PPM simulations and the result from the XFEM simulation in [61]. The results in Figure 8 show that the loading curves from the simulations in this study are consistent with the results from the XFEM method. The loading curve for the simulation with $\delta$ = 11.3 mm deviates from the results with $\delta$ = 7.5 mm and $\delta$ = 6 mm. Recall that the micromechanical modulus $C$ in (7) is an inverse function of the horizon. The larger peak load of the simulation using $\delta$ = 11.3 could be related to $C$ in the stabilization term. Figures 9 and 10 plot the contours of damage variable $\varphi$ and water pressure at $u_y$ = 2.5 ×$10^{-2}$ mm, respectively. The results in Figure 9 show that $\delta$ could slightly affect the damage zone and crack propagation length. Specifically, the lengths of the new crack are 93 cm, 97 cm, and 98 cm for the simulations with $\delta$ = 11.3 mm, 7.5 mm, and 6.0 mm, respectively. Figure 10 shows that the water pressure is negative (suction) around the propagated crack for all three cases. Consistent with the results in Figure 9, the value of $\delta$ affects the negative water pressure distribution around the newly formulated crack. The negative water pressure near the new crack can be caused by the competing factors between the crack opening space and the amount of water flowing into the fracture space. Following this reasoning, the negative water pressure could be generated by the crack opening and dilation in the bulk around the crack.

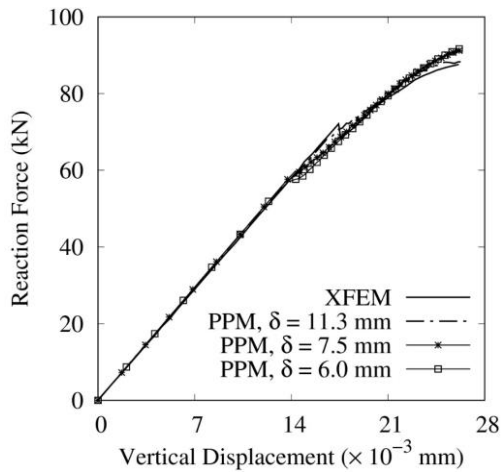

Figure 8: Comparison of the loading curves the simulations using three values of $\delta$ with the XFEM result [61].

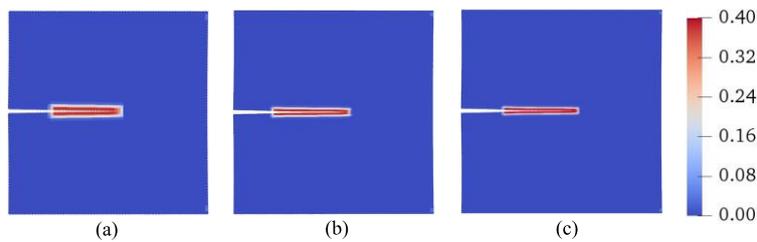

(a) (b) (c)

Figure 9: Contours of damage variable $\phi$ from the simulations with (a) $\delta$ = 11.3 mm, (b) $\delta$ = 7.5 mm, and (c) $\delta$ = 6 mm, at $u_y$ = 2.5 ×10$^{-2}$ mm (×50).

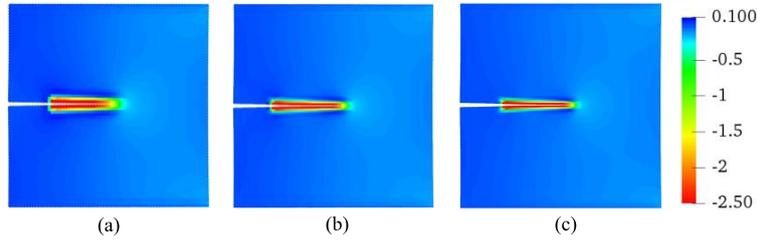

Figure 10: Contours of water pressure (MPa) from the simulations with (a) $\delta$ = 11.3 mm, (b) $\delta$ = 7.5 mm, and (c) $\delta$ = 6.0 mm, at $u_y$ = 2.5 ×10$^{-2}$ mm (×50).

### 4.1.2. Influence of m

The ratio $m = \delta/d$ is a measure of the number of mixed points within the horizon of a mixed point. It is another factor that could affect the numerical results with peridynamics. To study the influence of $m$ on the results, we rerun the simulations with $m$ = 2, 3, and 4 and $\delta$ = 7.5 mm. The spatial discretizations consist of 9,000, 20,000, and 36,000 mixed points, respectively. All other parameters remain the same.

Figure 11 plots the loading curves of the three simulations and the XFEM result [61]. Figures 12 and 13 compare the contours of damage and water pressure, respectively, at $u_y$ = 2.5 ×10$^{-2}$ mm.

In the early stage the loading curves in Figure 11 are identical for all simulations. However, the simulations using $m$ = 3 and $m$ = 4 predict larger peak loads than the XFEM result. As shown in Figure 12 (a)-(c) there is a noticeable decrease in crack length from simulations using larger values of $m$. For instance, the crack length is 135 mm from $m$ = 2 and the crack length is 90 mm from $m$ = 4. These differences in crack length are reflected in the force plot (Figure 11), e.g., a larger peak load for a larger value of $m$. It is noted that previous studies conducted using PPM [14] or hybrid FEM-PD (peridynamics) [34, 45] have successfully used $m$ = 2 to reproduce static and dynamic coupling phenomena in continuum analysis of porous media. From this example it appears that $m$ = 3 is an appropriate choice when modeling mode I crack propagation in unsaturated porous media.

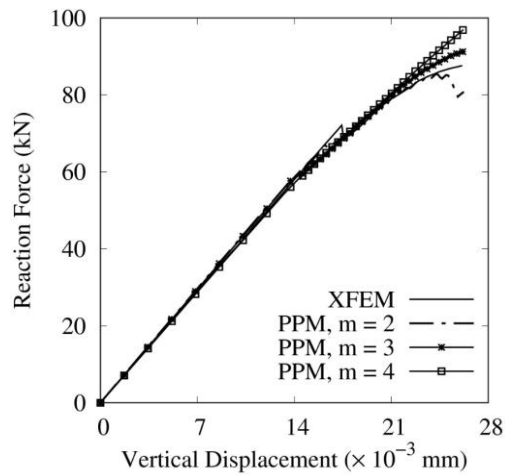

Figure 11: Comparison of loading curves from PPM simulations using a constant $\delta$ = 7.5 mm and three $m$-ratios to XFEM [61] results.

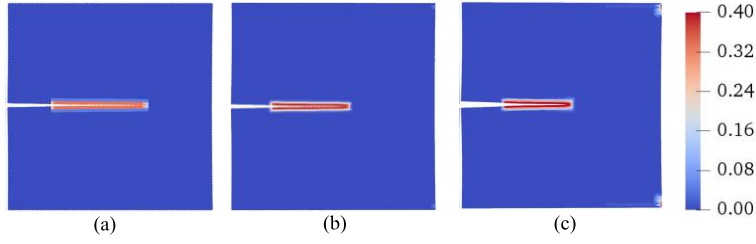

Figure 12: Contours of damage variable $\phi$ from the simulations with (a) $m = 2$, (b) $m = 3$, and (c) $m = 4$ at $u_y$ = 2.5 ×10$^{-2}$ mm (×50).

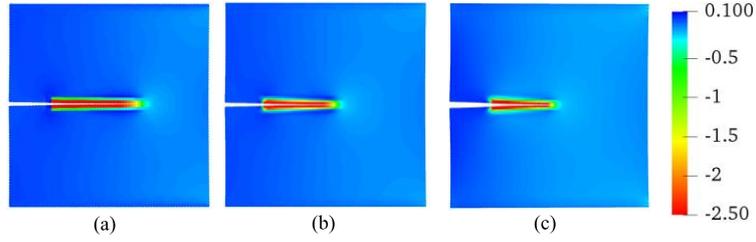

Figure 13: Contours of water pressure (MPa) from the simulations with (a) $m = 2$, (b) $m = 3$, and (c) $m = 4$ at $u_y$ = 2.5 ×10$^{-2}$ mm (×50).

## 4.2. Example 2: Mode II crack propagation

In this example we model mode II crack propagation by simulating shear loading of a pre-cracked porous body as inspired by an example in [6]. We note that the example in [6] was simulated by a saturated poromechanics model. Figure 14 depicts the problem domain, boundary conditions, and loading protocol. All fluid phase boundaries are assumed impermeable. The porous body is discretized into 19,200 mixed points with $d$ = 0.5 mm. All material parameters are the same as adopted in example 1. The specimen is prescribed a zero initial effective stress and water pressure as assumed in [6]. As shown in Figure 14 the loading rate $\dot{u}_x$ = 1 ×10$^{-3}$ mm/s. The total loading time $t$ = 30 s and $\Delta t$ = 0.02 s. We first present the base simulation results to show the crack propagation and water pressure variation in the specimen under the shear loading. Figure 15 presents the snapshots of the contour of damage variable $\phi$ superimposed on the deformation configuration at three loading stages. Figure 16 plots the contours of water pressure at three load stages. The results in Figure 15 indicate that the crack under the shear loading propagates upward following a slightly curved path. As shown in Figure 16 the water pressure is negative (i.e., matric suction) around the newly propagated crack and the area above the initial crack. Next, we repeat the base simulation with different horizons and $m$ to study the influence of $\delta$ and $m$ on the numerical results.

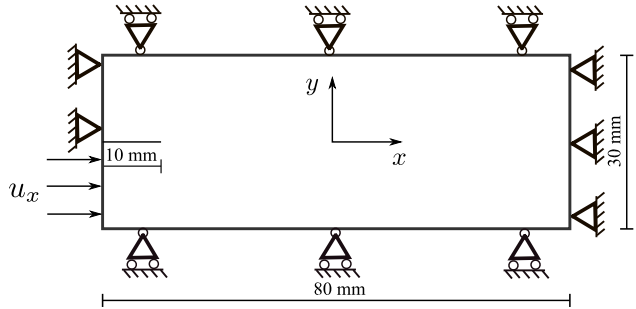

Figure 14: Problem setup for example 2.

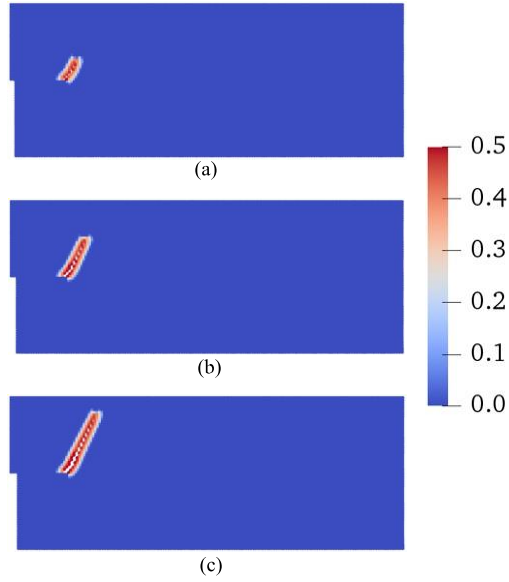

Figure 15: Contours of damage variable $\phi$ at $u_x$ = (a) 2.0 ×10⁻² mm, (b) 2.5 ×10⁻² mm, and (c) 3.0 ×10⁻² mm (× 50).

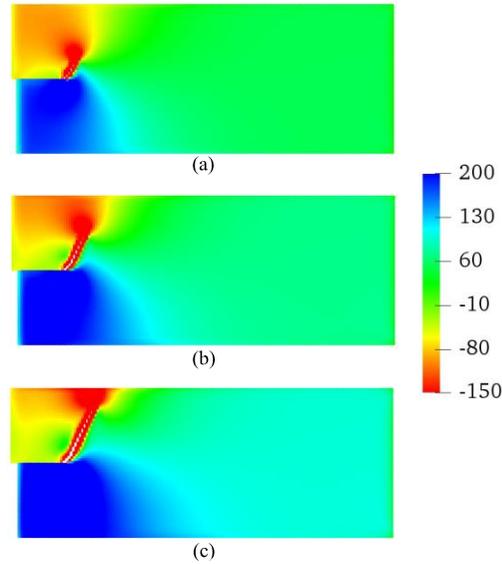

Figure 16: Contours of water pressure (kPa) at $u_x$ (a) 2.0 ×10⁻² mm, (b) 2.5 ×10⁻² mm, and (c) 3.0 ×10⁻² mm (× 50).

### 4.2.1. Influence of $\delta$

To study the influence of $\delta$ on the model II crack propagation, we conduct the simulations using three different horizons $\delta$ = 1.5 mm, 1.2 mm and 0.9 mm respectively with $m$ = 3. The discretizations consist of 19,200, 31,000, and 53,000 mixed material points, respectively. The results are shown in Figures 17, 18, and 19.

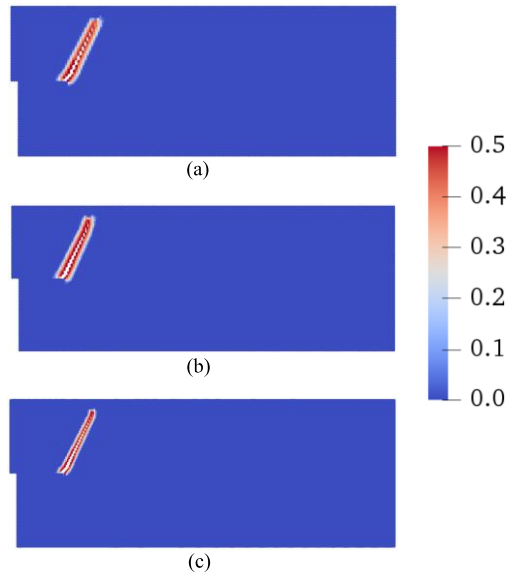

Figure 17: Contours of damage variable $\varphi$ from the simulations using (a) $\delta$ = 1.5 mm, (b) $\delta$ = 1.2 mm, and (c) $\delta$ = 0.9 mm at $u_x$ = 3 × 10$^{-2}$ mm (× 50).

Figure 17 presents the contours of damage variable $\varphi$ for the three simulations superimposed on the deformed configuration at $u_x$ = 6 ×10$^{-2}$ mm. The results in Figure 17 show that the crack propagation follows a similar path for the three simulations. This is consistent with the results in Figure 19 that compares the crack path for the three cases. Figure 18 plots the contours of water pressure in the problem domain for the three simulations. For all three cases the water pressure around the crack and above the initial crack is negative (suction) and the water pressure below the initial crack is positive.

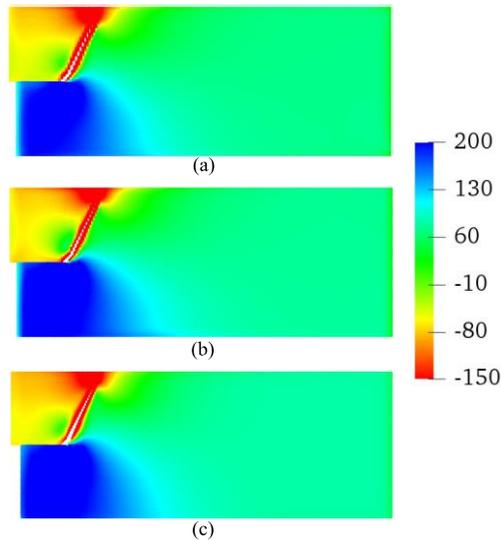

Figure 18: Contours of water pressure (kPa) from the simulations using (a) $\delta$ = 1.5 mm, (b) $\delta$ = 1.2 mm, and (c) $\delta$ = 0.9 mm at $u_x$ = 3 × 10$^{-2}$ mm (× 50).

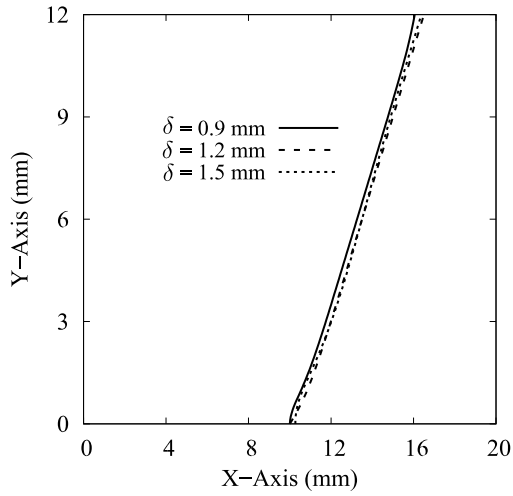

Figure 19: Comparison of the crack path from simulations using different $\delta$ and $m = 3$.

*4.2.2. Influence of m*

To study the influence of $m$ on the model II crack propagation in porous media we conduct simulations using $m$ = 3, 4, and 5 with $\delta$ = 1.5 mm. The discretizations consist of 19,200, 31,000, and 53,000 mixed material points, respectively. Figure 20 presents the contours of damage for the simulations at $u_x = 3 \times 10^{-2}$ mm. Figure 21 plots the contours of water pressure for the simulations. The crack paths from the three cases are compared in Figure 22.

It is indicated from Figure 22 that the choice of $m$ has a mild impact on the crack propagation path. However, the value of $m$ has a noticeable impact on both the water pressure in the problem domain as shown in Figure 21. Comparison of the results in Figures 19 and 22 shows that $m = 3$ could be a sufficient choice for modeling unguided mode II cracks in variably saturated porous media. We note that more studies should be conducted to guide the appropriate choice of $m$ for realistically modeling arbitrary cracks in unsaturated porous media through the proposed framework.

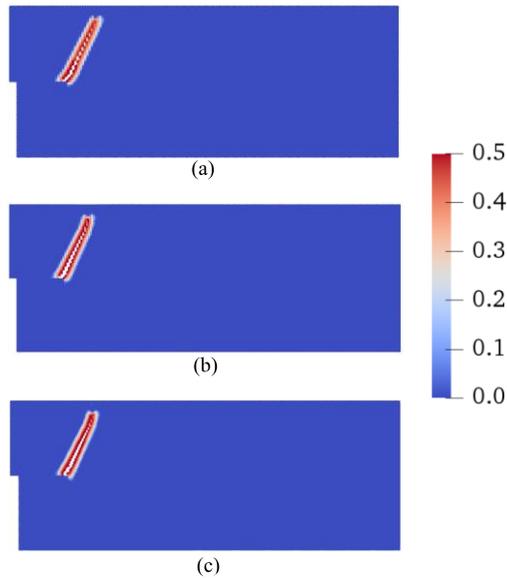

Figure 20: Contours of damage variable $\varphi$ from the simulations using (a) $m = 3$, (b) $m = 4$, and (c) $m = 5$ at $u_x = 3 \times 10^{-2}$ mm (×50).

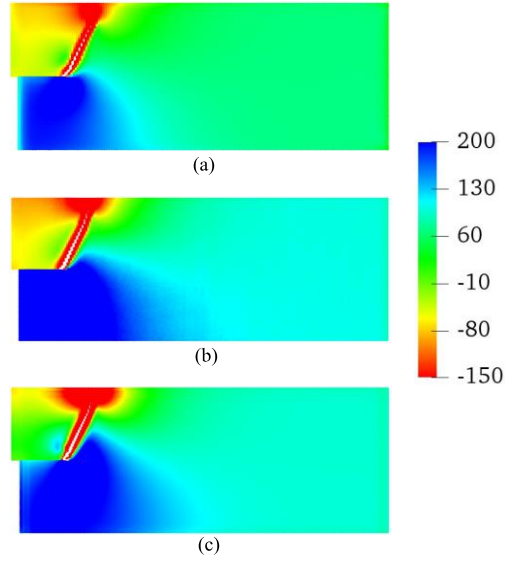

Figure 21: Contours of water pressure (kPa) from the simulations using (a) *m* = 3, (b) *m* = 4, and (c) *m* = 5 at $u_x$ = 3 ×10$^{-2}$ mm (×50).

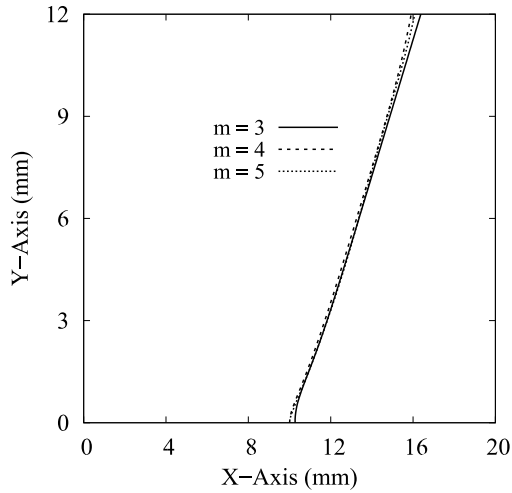

Figure 22: Comparison of the crack path from simulations using (a) *m* = 3, (b) *m* = 4, and (c) *m* = 5 at $u_x$ = 3 ×10$^{-2}$ mm.

### 4.3. Example 3: Wing crack propagation

In this example we simulate the wing crack propagation in an unsaturated porous medium through the implemented unsaturated fracturing periporomechanics model. Figure 23 depicts the problem geometry and loading protocol. The problem domain is discretized into 20,000 uniform mixed material points with *d* = 1 mm. The material parameters adopted are *K* = 70×10$^4$ kPa, $\mu_s$ = 15×10$^4$ kPa, $\phi_0$ = 0.33, $k_w$ = 1×10$^{-15}$ m$^2$, *n* = 1.5, $s_a$ = 500 kPa. The horizon $\delta$ = 3.05 *d*. The porous body is prescribed an initial uniform effective stress -49.5 kPa and initial suction 50 kPa (i.e., the zero initial total stress) $S_r$ = 0.99 from (16). The load rate $\dot{u}$ = 1.41 × 10$^{-5}$ mm/s as shown in Figure 23. The total loading time *t* = 3000 s and $\Delta t$ = 1 s.

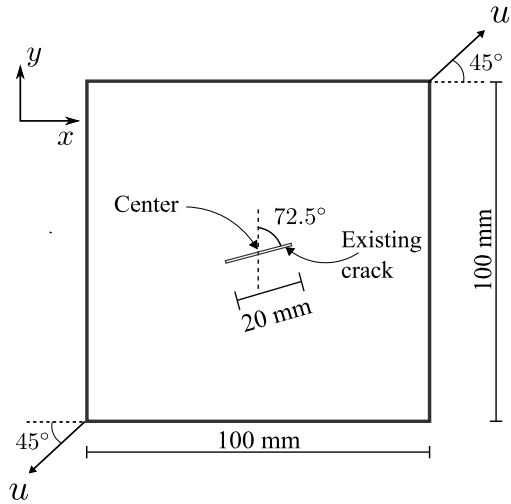

Figure 23: Problem setup for example 3.

Figure 24 plots the loading curve versus the applied horizon displacement. The vertical dash lines denote the applied displacements $u_x = u_y$ (a) $1.0 \times 10^{-2}$ mm (b) $1.1 \times 10^{-2}$ mm and (c) $1.2 \times 10^{-2}$ mm. The snapshots of contours of the damage variable $\varphi$ and water pressure at these three loading stages are presented in Figures 25 and 26, respectively. The results in Figure 25 show that the propagation path of the wing cracks is parallel to the diagonal line of the specimen. This crack propagation path matches the experimental result and the wing crack propagation in a single-phase solid under the same loading condition through peridynamics (e.g., [43, 46]). The contours of water pressure in Figure 26 indicate that the crack propagation has led to the increase of matric suction around the wing cracks. This may imply that the formation rate of new fracture space is larger than the rate of the volume of water flowing into the newly formed fracture space. Therefore, the degree of saturation in the fracture space is smaller than that in the bulk (see (16)). The overall low hydraulic conductivity (i.e., intrinsic permeability multiplied by the relative permeability) in the bulk under unsaturated condition could be a contributing factor for the low flow rate of water from the bulk into the fracture space. In what follows, we study the influence of initial matric suction and intrinsic permeability on the wing crack propagation by repeating the simulation with different initial matric suctions and intrinsic permeabilities.

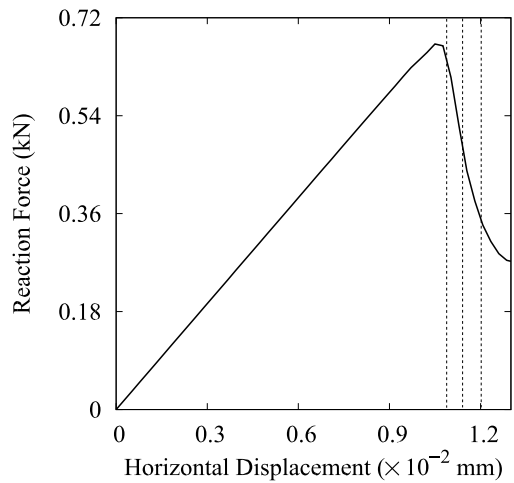

Figure 24: Plot of the loading curve over applied horizontal displacement.

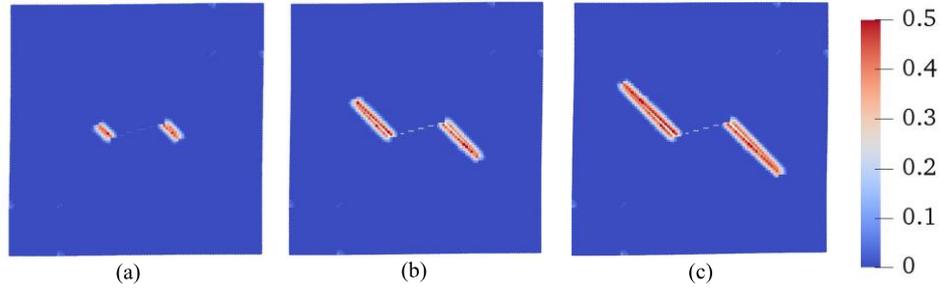

Figure 25: Contours of damage variable $\varphi$ superimposed on the deformed configuration at $u_x$ = (a) 1.0 ×10$^{-2}$ mm, (b) 1.1 ×10$^{-2}$ mm, and (c) 1.2 ×10$^{-2}$ mm (×50).

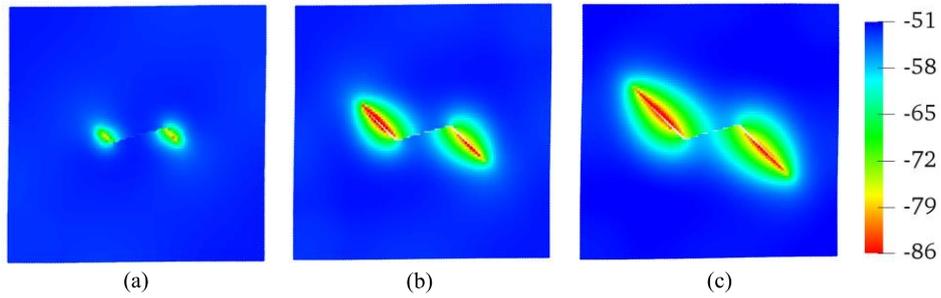

Figure 26: Contours of water pressure (kPa) superimposed on the deformed configuration at $u_x$ = (a) 1.0 ×10$^{-2}$ mm, (b) 1.1 ×10$^{-2}$ mm, and (c) 1.2 ×10$^{-2}$ mm (×50).

### 4.3.1. Influence of initial matric suction

To study the influence of initial matric suction on the wing crack propagation, the simulations were repeated with another two values of initial matric suction $s_2 = 100$ kPa and $s_3 = 200$ kPa. To maintain the zero initial total stress condition the corresponding initial effective stress is -96 kPa and -180 kPa, respectively. The degrees of saturation $S_r$ are respectively 0.96 and 0.9. All other parameters remain unchanged.

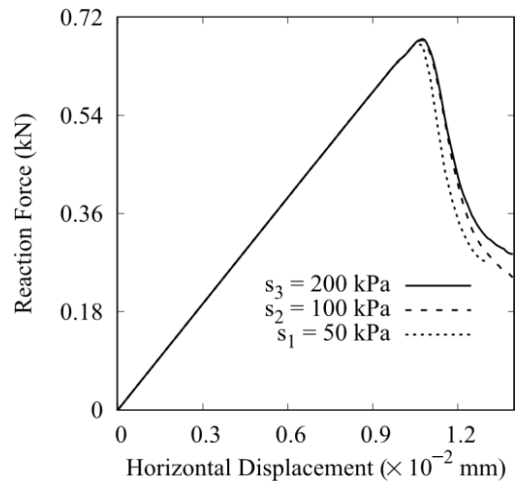

Figure 27: Comparison of the loading curve from the simulations with different initial suctions

Figure 27 plots the loading curves from the simulations with the three initial suctions. The initial matric suction has mild impact on the peak loading. Once the crack propagates, the reaction load for the specimen with larger matric suction is generally larger than that with smaller initial matric suction over the loading process. This could be due to the larger tensile strength of the specimen with larger matric suction. Figures 28 and 29 present the contours of damage and water pressure superimposed on the deformed configuration at the same loading stage. The comparison in Figure 28 show that the length of the wing cracks are associated with the initial suction under the same mechanical loading stage. The wing cracks are shorter in the specimen with larger initial matric suction. The magnitude of the initial suction does not impact the direction of the wing cracks due to the isotropic assumption. As shown in Figure 29 matric suction increases around the wing cracks for the simulations with different initial matric suctions.

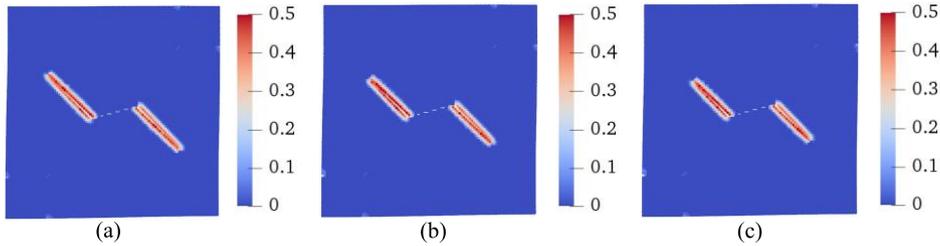

Figure 28: Contours of damage variable $\varphi$ from the simulations with the initial suction (a) $s_1$ = 50 kPa, (b) $s_2$ = 100 kPa, and (c) $s_3$ = 200 kPa at $u_x$ = 1.2 × 10$^{-2}$ mm (×50).

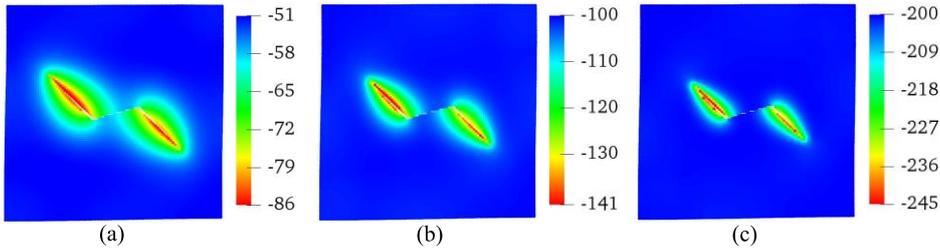

Figure 29: Contours of water pressure (kPa) from simulations with initial suction (a) $s_1$ = 50 kPa, (b) $s_2$ = 100 kPa, and (c) $s_3$ = 200 kPa at $u_x$ = 1.2 × 10$^{-2}$ mm (×50).

*4.3.2. Influence of intrinsic permeability*

We study the influence of intrinsic permeability on the wing crack propagation in unsaturated porous media by comparing the results from the simulations with three intrinsic permeabilities $k_1$ = 1 × 10$^{-13}$ m$^2$, $k_2$ = 1 × 10$^{-14}$ m$^2$, and $k_3$ = 1 × 10$^{-15}$ m$^2$. Figure 30 plots the loading curves from the three simulations. The loading curves are almost identical before the peak value of the simulation with $k_1$ = 1 × 10$^{-13}$ m$^2$. After that the loading curves for the three simulations diverge from each other. For the two simulations with smaller intrinsic permeabilities the specimen has a larger peak load. For all three cases, the crack starts propagating at the peak load. Figures 31 and 32 plot the snapshots of the contours of damage variable and water pressure at different loading stages, respectively. The results show that intrinsic permeability could have a dramatic influence on the length of wing cracks, as well as the matric suction around the wing cracks. As shown in Figure 31, at the same loading stage the length of the wing cracks from the simulation with smaller intrinsic permeability is smaller than the simulation with larger intrinsic permeability. The larger matric suction has generated around the wing crack for the simulation with smaller intrinsic permeability as shown in Figure 32. The larger matric suction may explain the reduced wing crack propagation for the simulation with smaller intrinsic permeability. The results have demonstrated that the intrinsic permeability could have a significant impact on the wing crack propagation and the matric suction around the wing crack in unsaturated porous media.

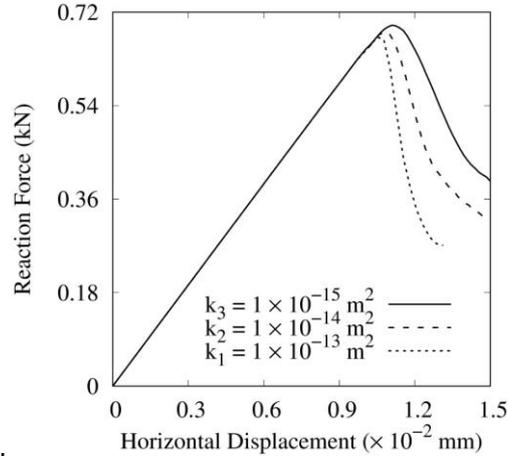

Figure 30: Comparison of the loading curves for the simulations with three intrinsic permeabilities.

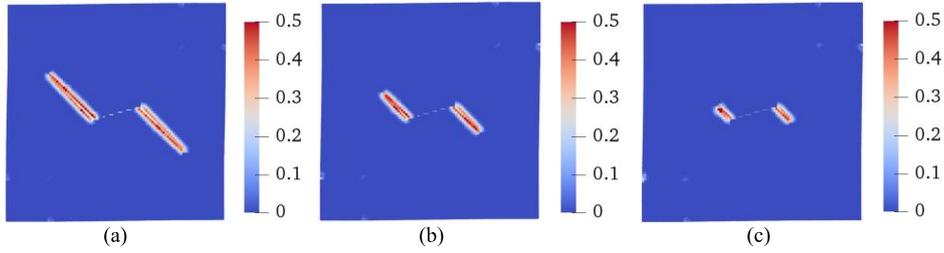

Figure 31: Contours of the damage variable $\varphi$ for the simulations with three intrinsic permeabilities (a) $k_1 = 1 \times 10^{-13}\,\text{m}^2$, (b) $k_2 = 1 \times 10^{-14}\,\text{m}^2$, and (c) $k_3 = 1 \times 10^{-15}\,\text{m}^2$ in the deformed configuration at $u_x = 1.2 \times 10^{-2}\,\text{mm}$ (×50).

### 4.4. Example 4: Non-planar cracking in an unsaturated parallelepiped soil specimen

In this example we simulate the formation of non-planar cracks in an unsaturated parallelepiped soil specimen triggered by matric suction variations. Figure 33 sketches the problem setup and boundary conditions.

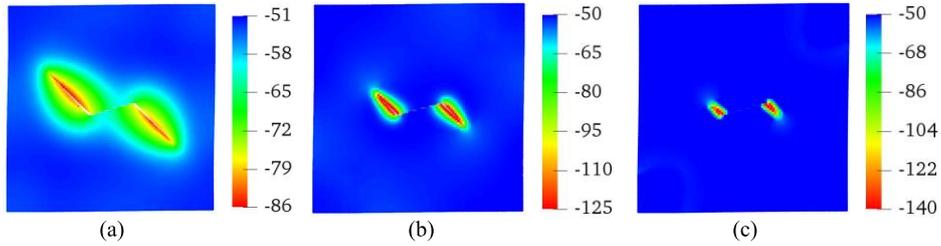

Figure 32: Contours of water pressure for the simulations with three intrinsic permeabilities (a) $k_1 = 1 \times 10^{-13}\,\text{m}^2$, (b) $k_2 = 1 \times 10^{-14}\,\text{m}^2$, and (c) $k_3 = 1 \times 10^{-15}\,\text{m}^2$ in the deformed configuration at $u_x = 1.2 \times 10^{-2}\,\text{mm}$ (×50).

The problem domain is discretized into 20,000 mixed points with $d = 1.8$ mm. The material parameters adopted are $K = 5 \times 10^3$ kPa, $\mu_s = 2.1 \times 10^3$ kPa, $\phi_0 = 0.1$, $G_c = 0.5$ J/m$^2$, $k_w = 1 \times 10^{-14}$ m$^2$, $n = 1.25$, $s_a = 100$ kPa. The horizon $\delta = 3.05\,d$. The initial uniform effective stress in the specimen is assumed -47 kPa. The initial matric suction is assumed 50 kPa and $S_r = 0.94$ from (16). In this case, the initial total stress in the specimen is null. A water flow boundary with a constant rate $q_y$ = 1.5 kg/(m³s) is imposed on the top surface through a fictitious boundary layer with the thickness of $\delta$. All other fluid boundaries are impermeable. The total simulation time $t$ = 3 hours (h) and $\Delta t$ = 1 s.

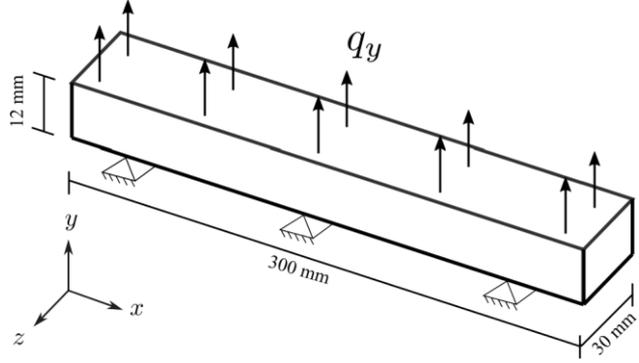

Figure 33: Problem setup of example 4.

First, we present the energy dissipation characteristics due to the formation of non-planar cracks in the specimen. We note the dissipated energy due to bond breakage is determined from (36) based on the effective force state concept. The dissipated energy density represents the energy consumed due to bond breakage at a material point. Figure 34 shows the total energy dissipation over the simulation time due to bond breakage. Here the total energy dissipation in the specimen is defined as the summation of the energy dissipation at all fracture points. Figure 35 presents the snapshots of the contour of dissipation energy associated with the bond breakage and the formation of non-planar cracks at three load steps denoted in Figure 34. The results demonstrate that the rate of energy dissipation due to bond breakage has increased dramatically upon the inception and propagation of cracks.

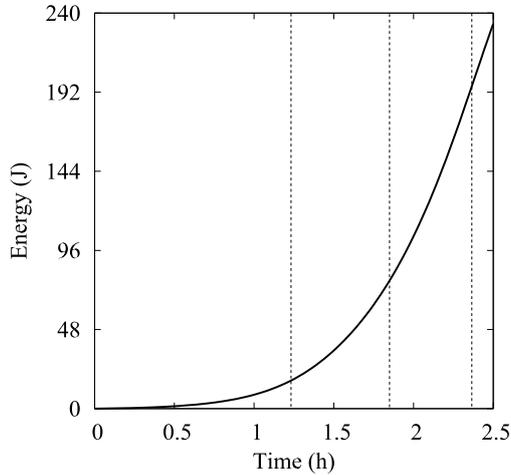

Figure 34: Plot of the total energy dissipation over time due to the crack formation.

Second we present the zoom-in snapshots of non-planar crack topology and the contour of matric suction in the specimen. The results are shown in Figures 36, 38, and 37. Figure 36 plots the contours of the horizontal displacement and the zoom-in crack topology at three time steps. Figure 37 depicts the snapshots of the dissipation energy density on the crack surface in the deformed configuration. Two major non-planar cracks have formed in the specimen due to the shrinkage in the *x* direction. Figure 36 demonstrates that the two major cracks initiate on the top and side surfaces and then propagate downward within the specimen. Figure 38 plots the contours of water pressure or matric suction on the crack surface. The matric suction on the top surface is the largest due to the outward water flow boundary condition. From the results in Figures 36 and 38 we could conclude in this example the non-planar cracks within the specimen has been triggered.

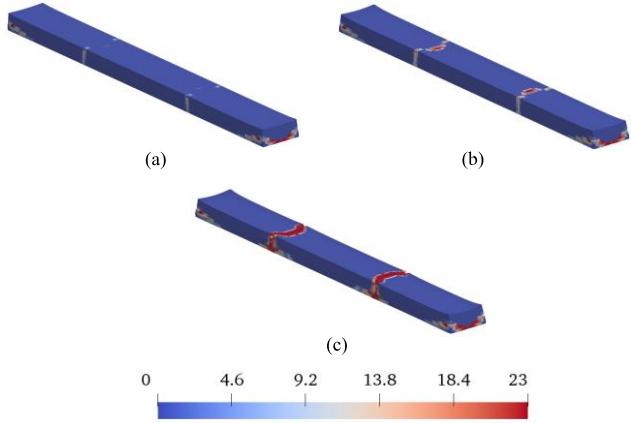

Figure 35: Contours of the energy dissipation density due to bond breakage (MJ/m$^3$) superimposed on the deformed configuration at (a) $t$ = 1.2 h, (b) $t$ = 1.8 h, and (c) $t$ = 2.3 h.

by the increasing matric suction from the top surface. Through this example, we have demonstrated the efficacy of the proposed fracturing periporomechanics framework for modeling unguided non-planar cracking in unsaturated porous media.

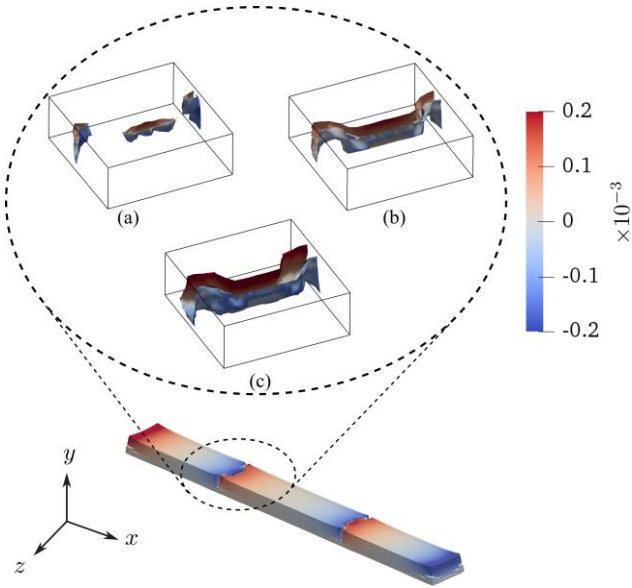

Figure 36: Snapshots of the contour of the horizontal displacement $u_x$ (m) and zoom-in crack topology at (a) $t$ = 1.2 h , (b) $t$ = 1.8 h, and (c) $t$ = 2.3 h.

## 5. Closure

In this article, as a new contribution we have formulated and implemented an unsaturated fracturing periporomechanics framework for unguided cracking in unsaturated porous media. In this new coupled periporomechanics model crack nucleation and propagation are modeled by an effective force based criterion that incorporates the effect of suction on cracking. It is noted that crack formation is completely autonomous and requires no external criteria nor a priori knowledge/assumption of the crack path. Unsaturated fluid flow in the fracture space is modeled via a simplified formation in line with the unsaturated fluid flow in bulk. A fractional step algorithm via the celebrated two-stage operator split and two-phase mixed meshless method are utilized to solve the coupled fracturing periporomechanics model. At each time step, an undrained deformation stage of skeleton by fixing fluid flow is solved first and then unsaturated fluid mass transport is solved in the updated fractured skeleton configuration. The implementation of the coupled fracturing periporomechanics model validated using numerical examples based on the extended finite element method in the literature. Through the numerical simulations, we have demonstrated the capability and robustness of the proposed nonlocal poromechanics in modeling Mode I and Mode II crack propagation in unsaturated porous media. The efficacy of the fracturing unsaturated periporomechanics model for non-planar cracks has been demonstrated by modeling unguided cracking in a three-dimensional soil specimen triggered by matric suction variations.

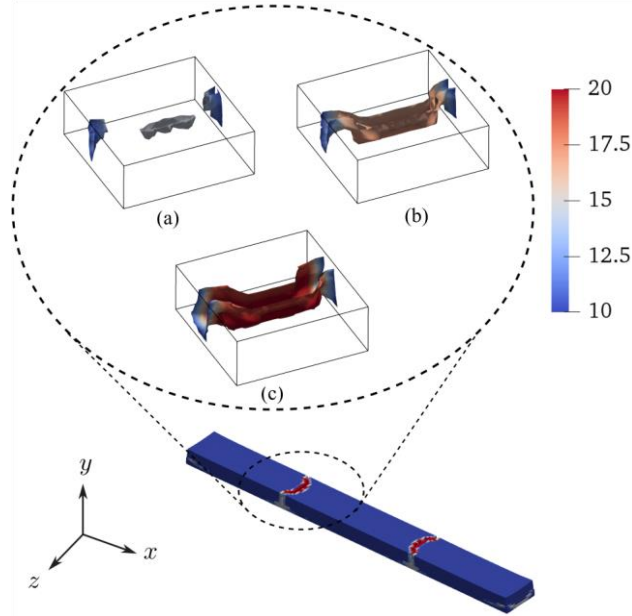

Figure 37: Snapshots of the contour of dissipation energy density (MJ/m$^3$) superimposed on the zoom-in crack topology at (a) $t$ = 1.2 h, (b) $t$ = 1.8 h, and (c) $t$ = 2.3 h.

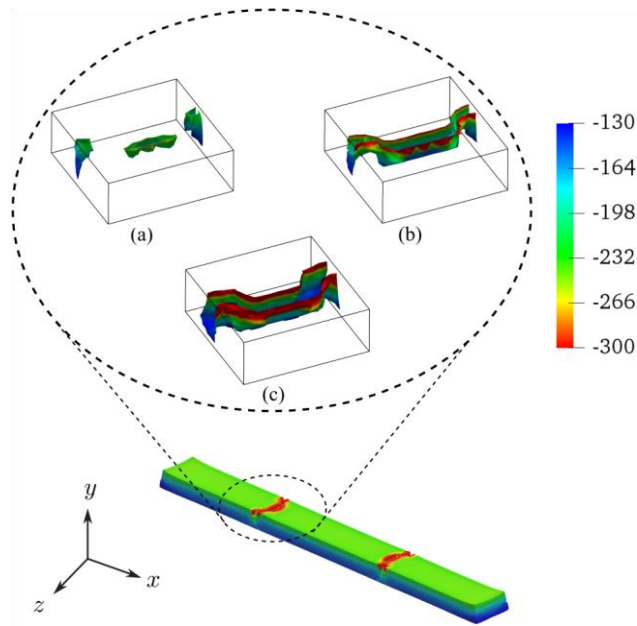

Figure 38: Snapshots of the contour of water pressure (kPa) superimposed on the zoom-in crack topology in the deformed configuration at (a) $t$ = 1.2 h , (b) $t$ = 1.8 h, and (c) $t$ = 2.3 h.

## Acknowledgments

This work has been supported by the US National Science Foundation under contract numbers 1659932 and 1944009.